\input amstex
\documentstyle{amsppt}
\magnification 1200
\NoRunningHeads
\NoBlackBoxes
\document

\def\Sym{\text{Sym}}

\def\tlo{\hat\otimes}

\def\tio{\tilde\otimes}
\def\ho{\hat\otimes}
\def\ad{\text{ad}}
\def\Hom{\text{Hom}}

\def\a{\frak a}

\def\Ua{U_q(\tilde\g)}
\def\U2{{\Ua}_2}
\def\g{\frak g}

\def\C{\Bbb C}

\def\<{\langle}
\def\>{\rangle}
\def\o{\otimes}
\def\e{\varepsilon}

\def\End{\text{End}}

\topmatter
\title Quantization of Lie bialgebras, I
\endtitle
\author {\rm {\bf Pavel Etingof and David Kazhdan} \linebreak
\vskip .1in
Department of Mathematics\linebreak
Harvard University\linebreak 
Cambridge, MA 02138, USA\linebreak
e-mail: etingof\@math.harvard.edu\linebreak kazhdan\@math.harvard.edu}
\endauthor
\endtopmatter

\vskip .05in
\centerline{\bf Abstract}         
\vskip .05in

In the paper  \cite{Dr3}
 V.Drinfeld formulated a number of problems in quantum group theory.
In particular, he raised the question about the existence of a 
universal quantization for Lie bialgebras, which arose from the
problem of quantization of Poisson Lie groups. When the paper 
\cite{KL} appeared Drinfeld
asked whether the methods of \cite{KL} could be useful for the problem of
universal quantization of Lie bialgebras. This paper gives a
positive  answer to a number of
Drinfeld's questions, using the methods and ideas of \cite{KL}. 
In particular, we show the existence of a universal
quantization. We plan to provide positive answers to  
most of the remaining questions in \cite{Dr3} in the following papers
of this series.

\vskip .05in
\centerline{\bf Introduction}
\vskip .05in

 The main result of this
paper is a construction of a
universal quantization for Lie bialgebras(see \cite{Dr3} Section 1). 

The paper consists of two parts. In the first part we construct the
quantization of a finite dimensional Lie bialgebra. 
In the second part we generalize this result to the
infinite-dimensional case. The 
construction in the first part consists of three steps. 

1) Given a finite dimensional Lie bialgebra $\a$ over a 
field $k$ of characteristic zero, we construct the 
double $\g$  of $\a$. Our definition of the double 
coincides with the one in  \cite{Dr1}. 
We consider the category $\Cal M$
whose objects are $\g$-modules and
$\text{Hom}_{\Cal M}(U,W)=\text{Hom}_\g(U,W)[[h]]$. For any associator
$\Phi$ (\cite{Dr2,Dr4}) we define
a structure of a braided monoidal category on $\Cal M$, as in \cite{Dr2}.

2) We construct Verma modules $M_+$, $M_-$ over $\g$,
and use them to construct a fiber functor from $\Cal M$ 
to the tensor category of topologically free $k[[h]]$
modules: $F(V)=\Hom_{\Cal M}(M_+\o M_-,V)$.
 According to the categorical yoga, the existence of
such a functor implies the existence of 
 a (topological) Hopf algebra $H$ isomorphic to $U(\g)[[h]]$
such that the tensor category $\Cal M$ is
equivalent to the category  of representations of $H$. We show that $H$ is
isomorphic, as a topological algebra,
 to $U(\g)[[h]]$,
where $U(\g)$ is the universal enveloping algebra of the Lie algebra
$\g$. 

3) We construct Hopf subalgebras $H_\pm$ of $H$ and
show that $H_+$ is a quantization of $\a$ and
that the algebra $H$ is the quantum double of the Hopf algebra $H_+$. 

{\bf Remark.} We do not expect the existence of a quantization
of any Lie bialgebra $\a$ which is isomorphic to $U(\a)[[h]]$ as a topological 
algebra.

As an application of our techniques, 
we prove that any classical r-matrix $r$
over an associative algebra $A$ ($r\in A\o A$)
can be quantized. In other words, there exists a 
quantum R-matrix $R\in A\o A[[h]]$ such that $R=1+hr$. We also show 
that $R$ is unitary ($R^{21}R=1$) if $r$ is unitary ($r^{21}=-r$).
This answers questions in Section 3 of \cite{Dr3}.
As another
application, we show the existence of the quantization of a quasitriangular 
Lie bialgebra $\a$ (not necessarily finite dimensional)
such that the obtained quantized universal enveloping
algebra has a quasitriangular structure and is isomorphic to $U(\a)[[h]]$ 
as a topological algebra, which solves questions in Section 4 of \cite{Dr3}.

The construction of quantization given in Part I has two drawbacks.
First, it does not work literally for infinite dimensional Lie bialgebras.
Second, it does not allow to prove functoriality and universality of 
quantization. Therefore, in Part II we slightly modify the construction,
which puts the results of the first part 
in a more general setting. Now we consider arbitrary Lie bialgebras, not 
necessarily finite dimensional. In this case the double $\g$ of $\a$
can also be constructed, but it carries a nontrivial topology if 
$\text{dim}\a=\infty$. Instead of the category of all $\g$-modules, we now 
consider the category $\Cal M^e$ whose objects are
 equicontinuous $\g$-modules, which are 
topological $\g$-modules satisfying certain conditions. On this category,
we define a braided monoidal structure analogously to the finite-dimensional 
case.

We construct Verma modules $M_+$, $M_-$ over $\g$ analogously to
the finite-dimensional case. The module $M_-$ is equicontinuous.
The module $M_+$, in general,
is not equicontinuous, but the module $M_+^*$, dual to $M_+$ in an appropriate
topology,
is an equicontinuous $\g$-module. Using $M_-$ and $M_+^*$, we define a fiber functor 
from $\Cal M^e$ to the category of topological $k[[h]]$-modules,
by $F(V)=\Hom_{\Cal M^e}(M_-,M_+^*\o V)$. Since the module $M_+$ is
not always equicontinuous, this functor is not always representable in $\Cal
M^e$. We define a tensor structure on
$F$ similarly to the finite dimensional case, and show that if $\g$ is finite
dimensional, the functors
obtained in the first and second parts of the paper are isomorphic as
tensor functors. 

Next, we consider the algebra $H=\End F$. It is a topological 
algebra over $k[[h]]$ with a ``coproduct'' $\Delta$, which maps $H$ into
 a completion of $H\o H$, but not necessarily in $H\o H$.  

Finally, we construct a subalgebra $H_+$ of $H$ such that 
$\Delta(H_+)\subset H_+\o H_+$. This is a quantized universal enveloping 
algebra which is a quantization of $\a$. For finite dimensional $\a$, this
quantization is isomorphic to the one obtained in the first part.

At the end of the paper we settle Drinfeld's question of 
the existence of a universal quantization of Lie bialgebras by 
showing that the quantization
obtained in the second part of the paper is universal.
In Drinfeld's language this means that the product
and coproduct in the quantized algebra express
in terms of acyclic tensor calculus via the commutator
and cocmmutator. This result implies
that our quantization of Lie bialgebras is a functor from the category
of Lie bialgebras to the category of topological Hopf algebras. It 
also shows that our quantizations of classical r-matrices, 
unitary r-matrices, and quansitriangular Lie bialgebras are universal
and functorial.
Thus we 
answer positively the corresponding questions of Drinfeld \cite{Dr3}. 

{\bf Remarks.} 
1. The material of Part I does not seem sufficient for proving the 
universality and functoriality. 
In fact, during the computation of the $h^2$-term of
multiplication in $U_h(\a)$, using the method of Part I, 
one gets non-acyclic expressions, which cancel at the end of computation.
Thus, the generalization to the infinite-dimensional case
is essential for the proof of functoriality, even for
finite-dimensional Lie
bialgebras. 

2. Most of the results of the paper could be formulated and proved over 
the ring
$k[h]/(h^N)$ rather than $k[[h]]$, and then the results over $k[[h]]$ could be 
obtained as a limit. The only problem arises with the notions
of the dual quantized universal enveloping algebra and the quantum double, 
which collapse
over $k[[h]]/(h^N)$. This is why we chose to work over $k[[h]]$. 

In fact, it is easy to see that the main results of the paper hold 
in a more general setting than stated. Namely, one can take the Lie bialgebra
$\a$ to be ``dependent on $h$'', i.e. to be a Lie bialgebra over the 
ring $k[[h]]$, which is topologically free as a $k[[h]]$-module. 
The universal acyclic formulas for quantization whose existence is
shown in Section 10 are well defined for this case, and define
a functor $\a\to U_h(\a)$, from the category of
Lie bialgebras over $k[[h]]$ which are topologically free 
as $k[[h]]$-modules to the category of quantized universal
enveloping algebras. In the second paper of this series 
we will show that this functor
is in fact an equivalence of categories.

Moreover, it is easy to see that the acyclic formulas of Chapter 10
which express the product and coproduct in $U_h(\a)$ in terms of 
commutator $[,]$ and cocommutator $\delta$ of $\a$, are in fact 
formal series in $[,]$ and $h\delta$ with coefficients 
independent of $h$. This allows to generalize the results even further.
Namely, let $K$ be any local Artinian or pro-Artinian 
algebra over $\Bbb Q$, and $I$ be the maximal ideal in $K$.
Let $k=K/I$ (it is a field of characteristic 0). Given a Lie 
bialgebra $\a$ over $K$ which is (topologically) free as a $K$-module
and cocommutative modulo $I$. Then the quantization functor of 
Chapter 10 is defined and assigns to $\a$ a quantized universal
enveloping algebra $U_{quant}(\a)$ over $K$. 
In the second paper we will show that this is an equivalence
of categories between the category of Lie bialgebras over $K$
which are topologically free as $K$-modules, and the category
of quantized universal enveloping algebras over $K$.  
In particular, if $K=k[[h]]$ then $U_{quant}(\a)=U_h(\tilde\a)$,
where $\tilde\a$ is $\a$ with the same commutator and cocommutator
$\delta_{\tilde\a}=h^{-1}\delta_\a$. 

The third paper of this series is not written yet. Therefore we will only
indicate the topics which we are planning to present in this part. First
of all, we plan to  consider the case of graded
bialgebras with finite-dimensional homogeneous components and to show that
in this case our formal quantization  defines a family of
Hopf algebras $H_h$, depending on a parameter $h\in k$. 
Our second goal is to
prove that for Kac-Moody bialgebras  our quantization
coincides with the quantum Kac-Moody algebra. As another application
of our techniques we plan to show how to 
define a quantum analog of the Kac-Moody algebra
for arbitrary symmetrizable Cartan
matrix (not necessarily integral)
and show that for generic values of $q$ the "size" of the
quantized algebra is the same as of the usual Kac-Moody algebra.
This would settle the questions in Section 8 of \cite{Dr3}.
 
\vskip .1in
\centerline{\bf Acknowledgements}
\vskip .1in

The authors are pleased to 
thank Maxim Vybornov for useful conversations and Dror
Bar-Natan for discussions and information about associators. 
The authors also want to thank the referee for many valuable suggestions.
The first author was partially supported by an NSF postdoctoral fellowship.

\vskip .1in
\centerline{\bf Part I} 
\vskip .1in
\centerline{\bf 1. Drinfeld category.} 
\vskip .1in

The definitions and statements of Sections 1.1, 1.2
can be found in \cite{Dr1}.

\vskip .05in
1.1. {\it Lie bialgebras.}

Throughout this paper, $k$ denotes a field of characteristic zero.
Let $\a$ be a Lie algebra over $k$, and. We $\delta$ be a linear map 
$\delta:\a\to\a\o\a$.

\proclaim{ Definition.} One says that the map $\delta$ defines 
a Lie bialgebra structure on $\a$ if it satisfies two conditions:

(i) $\delta$ is a 1-cocycle of $\a$ with coefficients in $\a\o\a$, i.e.  
$$
\delta([ab])=[1\o a+a\o 1,\delta(b)]+[\delta(a),1\o b+b\o 1];
$$

(ii) The map $\delta^*:\a^*\o \a^*\to\a^*$ dual to $\delta$ is a Lie
bracket on $\a^*$.

In this case $\delta$ is called the cocommutator of $\a$. 
\endproclaim

If $\a$ is a finite dimensional Lie bialgebra then $\a^*$
is a Lie bialgebra as well. Namely, the commutator in $\a^*$ is dual to the 
cocommutator in $\a$, and the cocommutator in $\a^*$ is dual to the commutator
in $\a$. If $\a$ is infinite-dimensional, then $\a^*$ is not in general a
Lie bialgebra but is a topological Lie bialgebra. That is, $\a^*$ is a Lie 
algebra in the usual sense, but the cocommutator
maps $\a^*$ into the completed tensor product $\a^*\hat\o \a^*$
and not necessarily into the usual tensor product $\a^*\o\a^*$. 

For any Lie bialgebra $\a$,
 the vector space $\g=\a\oplus\a^*$ has a natural structure of a Lie 
algebra. Namely, $\a,\a^*$ are Lie subalgebras in $\g$ with bracket defined 
above, and commutator between elements of $\a,\a^*$ is given by
$$
[a,b]=(\text{ad}^*a)b-(1\o b)(\delta(a)), a\in\a, b\in \a^*,\tag 1.1
$$ 
where $\text{ad}^*$ denotes the coadjoint action. There is an invariant
nondegenerate inner product on $\g$ given by $\<a+a',b+b'\>=a'(b)+b'(a)$, $a,b\in\a,a',
b'\in \a^*$. It is easy to show that (1.1) is the unique extension
of the commutator from $\a,\a^*$ to $\g$ for which the inner product $<,>$ is
ad-invariant.

\vskip .05in
1.2. {\it Manin triples.}

{\bf Definition.} A triple $(\g,\g_+,\g_-)$, where $\g$ is a finite dimensional
Lie algebra with a nondegenerate invariant
inner product $\<,\>$, and $\g_+$,$\g_-$ are isotropic Lie subalgebras,
such that $\g=\g_+\oplus\g_-$ as a vector space, 
is called a finite dimensional Manin triple. To every finite
dimensional Lie bialgebra $\a$ one can
associate the corresponding Manin triple $(\g=\a\oplus\a^*,\a,\a^*)$,
where the Lie structure on $\g$ is as above. Conversely, if $(\g,\g_+,\g_-)$
is a finite dimensional
 Manin triple then $\g_+$ (and $\g_-$) is naturally a Lie bialgebra.
Namely the pairing $\<,\>$ identifies $\g_+$ with $\g_-^*$, so we
can define $\delta:\g_+\to\g_+\o\g_+$ to be the 
dual map to the commutator of $\g_-$. 
This map is a 1-cocycle
of the Lie algebra $\g_+$ with coefficients in the module $\g_+\o\g_+$, 
so it defines a structure of a Lie bialgebra on $\g_+$. 

Thus, there is a one-to-one correspondence between finite dimensional
Lie bialgebras
and finite dimensional Manin triples. 

If $\a$ is a Lie bialgebra then the Lie algebra
 $\g=\a\oplus\a^*$ also has a natural structure
of a Lie
bialgebra. Namely, the cocommutator on $\g$ is 
$\delta_\g=\delta_\a\oplus(-\delta_{\a^*})$, where $\delta_\a$, $\delta_{\a^*}$
are the cocommutators of $\a$, $\a^*$. 

The 1-cocycle $\delta_\g$
is the coboundary of an element in $\g\o\g$.  
Namely, if $r\in \a\o\a^*\subset\g\o\g$ is the canonical
element corresponding to the identity operator $\a\to\a$, 
then $\delta_\g=dr$, where $r$ is regarded 
as a 0-cochain of $\g$ with coefficients in $\g\o\g$, and $d$ is the 
differential in the cochain complex;
that is $\delta_\g(x)=[x\o 1+1\o x,r]$. 

The Lie bialgebra $\g$ is called the double of $\a$.

Let $\a$ be a Lie algebra, and $r\in\a\o \a$.
The equation
$$
[r_{12},r_{13}]+[r_{12},r_{23}]+[r_{13},r_{23}]=0\tag 1.2
$$
in $U(\a)^{\o 3}$ is called the classical Yang-Baxter equation.
It is easy to check that the canonical element $r$ satisfies this equation.
  
{\bf Definition.}
We say that a Lie bialgebra $\a$ is quasitriangular if it is equipped with
an element $r\in\a\o\a$ satisfying the classical Yang-Baxter equation,
such that $\delta(a)=[a\o 1+1\o a,r]$ for any $a\in\a$ (i.e. $\delta$ is 
a coboundary of $r$). 

For example, the double $\g$ of any finite dimensional Lie bialgebra $\a$
equipped with the canonical element $r$ is a quasitriangular Lie bialgebra.   

\vskip .05in
1.3. {\it Associators.}
Recall some notation and definitions from the theory of associators
\cite{Dr2,BN}.  
Let $T_n$ be the algebra over $k$ generated by 
elements $t_{ij}$, $1\le i,j\le n$, $i\ne j$, with defining relations
$t_{ij}=t_{ji}$, $[t_{ij},t_{lm}]=0$ if $i,j,l,m$ are distinct, and
$[t_{ij},t_{ik}+t_{jk}]=0$. 

Let $P_1,...,P_n$ be disjoint subsets of $\{1,...,m\}$.
There exists a unique homomorphism 
$\rho_{P_1...P_n}: T_n\to T_m$ defined by
$$
\rho_{P_1...P_m}(t_{ij})=\sum_{p\in P_i,q\in P_j}t_{pq}.\tag 1.3
$$
For any $X\in T_n$, we denote $\rho_{P_1...P_n}(X)$
by $X_{P_1,...,P_n}$.

Let $\Phi\in T_3$. 
The relation 
$$
\Phi_{1,2,34}\Phi_{12,3,4}=\Phi_{2,3,4}\Phi_{1,23,4}\Phi_{1,2,3}\tag 1.4
$$
in $T_4[[h]]$ (=relation (1.2) in \cite{Dr2})
is called the pentagon relation.
 
Let $B=e^{ht_{12}/2}\in T_2[[h]]$. The relations
$$
\gather
B_{12,3}=
\Phi_{3,1,2}B_{1,3}\Phi_{1,3,2}^{-1}B_{2,3}\Phi_{1,2,3},\\
B_{1,23}=
\Phi_{2,3,1}^{-1}B_{1,3}\Phi_{2,1,3}B_{1,2}\Phi_{1,2,3}^{-1}.
\tag 1.5
\endgather
$$
in $T_3[[h]]$ (=relations (3.9a),(3.9b) in \cite{Dr2})
are called the hexagon relations. 

The element $\Phi$ is called an associator if it satisfies
the pentagon and hexagon relations. 

For $k=\C$,
an example of an associator is the Drinfeld associator $\Phi_{KZ}$
obtained from the KZ 
equations, as explained in \cite{Dr2}.  

The following theorem about associators is due to Drinfeld 
(\cite{Dr4}, Theorem A').

\proclaim{Theorem 1.1} There exists an associator defined over $\Bbb Q$. 
\endproclaim

This theorem implies that there exists an associator 
defined over any field $k$ of 
characteristic zero. From now on we will fix such an associator $\Phi$.

\vskip .05in
1.4. {\it Drinfeld category.}

Let $\g$ be a Lie algebra over $k$, and $\Omega\in S^2\g$ be 
a $\g$-invariant element. 

We will be mostly interested
in the case
when $\g$ belongs to a finite dimensional Manin triple $(\g,\g_+,\g_-)$,
and $\Omega=\sum_i g_i\o g^i$, where $\{g_i\}$ is a basis of $\g$,
and $\{\g^i\}$ is the dual basis to $\{\g_i\}$  
with respect to the invariant inner product on $\g$. 
In this case the element $\Omega$ is called the Casimir element.

Let $\Cal M$ denote the category whose objects are $\g$-modules, and
$\text{Hom}_{\Cal M}(U,W)=\text{Hom}_\g(U,W)[[h]]$.
This is a $k[[h]]$-linear additive category. For brevity we will later
write $\text{Hom}$ for $\text{Hom}_{\Cal M}$.
 
Drinfeld \cite{Dr2}
defined a structure of a braided monoidal category on $\Cal M$ as follows.

For any $V_1,V_2,V_3\in \Cal M$, consider a homomorphism
$\theta: T_3[[h]]\to
\End(V_1\o V_2\o V_3)$
by $\theta(t_{ij})=\Omega_{ij}$, and define
$\Phi_{V_1V_2V_3}=\theta(\Phi)$. 

For any $V_1,V_2\in \Cal M$, define
$V_1\o V_2\in\Cal M$ to be the usual tensor product of $V_1,V_2$
and the associativity morphism to be
$\Phi_{V_1V_2V_3}$, regarded as an element of
$\text{Hom}((V_1\o V_2)\o V_3,V_1\o (V_2\o V_3))$.
For any $V_1,V_2\in \Cal M$, introduce the braiding 
$\beta_{V_1V_2}:V_1\o V_2\to V_2\o V_1$
by the formula
$\beta=s\circ e^{h\Omega/2}$, where $s$ is the permutation. 
It follows from relations (1.4), (1.5)  
that the morphisms $\Phi_{V_1V_2V_3}$ and $\beta_{V_1V_2}$
define the structure of a braided monoidal category on $\Cal M$
(see \cite{Dr2}).

\vskip .1in
\centerline{\bf 2. The fiber functor.}
\vskip .1in

2.1. {\it The category of topologically free $k[[h]]$-modules.}

Let $V$ be a vector space over $k$. Then the space $V[[h]]$ 
of formal power series in $h$ with coefficients in $V$ has a natural
structure of a topological $k[[h]]$-module.
We call a topological $k[[h]]$-module topologically free if it is isomorphic 
to $V[[h]]$ for some $V$.

Let $\Cal A$ be the category of topologically free $k[[h]]$-modules,
where morphisms are continuous $k[[h]]$-linear maps. It is an 
additive category. Define the tensor structure on $\Cal A$ as follows:
for $V,W\in A$ define $V\o W$ to be the projective limit of
the $k[h]/h^n$-modules $(V/h^nV)\o_{k[h]/h^n}(W/h^nW)$ as $n\to\infty$. 
 
Let Vect be the category of vector spaces. We have the functor
of extension of scalars, $V\mapsto V[[h]]$, acting from 
$\text{Vect}$ to $\Cal A$. This functor respects 
the tensor product, i.e. $(V\o W)[[h]]$ is naturally
isomorphic to $V[[h]]\o W[[h]]$. 
The category $\Cal A$ equipped with the functor
$\o$ is a symmetric monoidal category.

If $X\in \Cal A$ then $X^*=\text{Hom}_{\Cal A}(X,k[[h]])$ is a 
topologically free $k[[h]]$-module. 
The assignment $X\to X^*$ is a contravariant functor from $\Cal A$
to itself.

\vskip .05in
2.2. {\it The forgetful functor.}

Let
$(\g,\g_+,\g_-)$ be a finite dimensional Manin triple, $\Omega\in S^2\g$ be 
the Casimir element associated to the inner product $\<,\>$ on $\g$, and
 $\Cal M$ 
be the Drinfeld category associated to $\g$.

Let $F:\Cal M\to \Cal A$ be the functor 
given by $F(M)=\text{Hom}(U(\g),M)$, where
$U(\g)$ is regarded as a left $\g$-module. 
This functor is naturally isomorphic to 
the ``forgetful'' functor which assigns to every $\g$-module $M$
the $k[[h]]$-module $M[[h]]$. The isomorphism between these
two functors is given by the assignment $f\in F(M)\to f(1)\in M[[h]]$.

\vskip .05in

2.3. {\it The Verma modules.}

Consider the Verma modules
$M_+=\text{Ind}_{\g_+}^\g \bold 1$, $M_-=\text{Ind}_{\g_-}^\g \bold 1$
(here $\bold 1$ denotes the trivial 1-dimensional representation).
By the Poincare-Birkhoff-Witt theorem,
the product in $U(\g)$ defines linear isomorphisms
$U(\g_+)\o U(\g_-)\to U(\g)$, and \linebreak $U(\g_-)\o U(\g_+)\to U(\g)$. 
This shows that the modules $M_{\pm}$ are freely generated over 
$U(\g_{\mp})$ by vectors $1_{\pm}$ such that $\g_{\pm}1_{\pm}=0$, 
and are identified (as vector spaces) with $U(\g_\mp)$
via $x1_{\pm}\to x$.

Since the vectors $1_\pm\o 1_\pm\in M_\pm\o M_\pm$ are $\g_\pm$-invariant,
there exist unique $\g$-module morphisms
 $i_{\pm}:M_{\pm}\to 
M_{\pm}\o M_{\pm}$ such that $i_\pm(1_\pm)=1_\pm\o 1_\pm$.
These morphisms in the category $\Cal M$ will play 
a crucial role in our constructions below.

\proclaim{Lemma 2.1} The assignment $1\to 1_+\o 1_-$
extends to an isomorphism of $\g$-modules $\phi: U(\g)\to M_+\o M_-$. 
\endproclaim

\demo{Proof} Since $M_\pm$ has been identified with $U(\g_\mp)$,
we can regard the map $\phi$ as a linear map $U(\g)\to U(\g_-)\o U(\g_+)$.
It is clear that this map preserves the standard filtration, so 
it defines a map of the associated graded objects: $S\g\to S\g_-\o S\g_+$.
This map is the isomorphism induced by the isomorphism $\g\to\g_-\oplus\g_+$.
Therefore, $\phi$ is an isomorphism. $\square$  
\enddemo

Lemma 2.1 implies that the functor $F$ can be identified with the functor 
$V\to \text{Hom}(M_+\o M_-,V)$. This 
definition of $F$ will be used from now on. 

\vskip .05in

2.4. {\it Tensor structure on the functor $F$.}

Let $(\Cal C,\o)$ be a monoidal category, 
$\Phi$ be the associativity constraint in $\Cal C$, and
$\bold 1$ be the identity object in $\Cal C$.
For simplicity we assume that 
$\bold 1\o X=X\o\bold 1=X$ for any object $X\in\Cal C$,
and the functorial isomorphisms $X\o X\o\bold 1$, $X\to \bold 1\o X$ the are 
identity morphisms.

Let $F:\Cal C\to\Cal A$
be a functor such that
$F(\bold 1)=k[[h]]$.

{\bf Definition.} 
By a tensor structure on the functor $F$ one means a functorial isomorphism
$J_{VW}:F(V)\o F(W)\to F(V\o W)$ satisfying
the associativity identity
$F(\Phi_{VWU})J_{V\o W,U}\circ (J_{VW}\o 1)=J_{V,W\o U}\circ (1\o J_{WU})$,
such that for any object $V$ $J_{V\bold 1}=J_{\bold 1V}=1$. 
A functor equipped with a tensor structure 
is called a tensor functor.
  
Now we describe a tensor structure on the functor
$F$ constructed in Section 2.2. 

For any $v\in F(V)$, $w\in F(W)$ define $J_{VW}(v\o w)$ to be the
composition of morphisms:
$$
\gather
M_+\o M_- @>i_+\o i_->> (M_+\o M_+)\o (M_-\o M_-) 
@>\text{associativity morphism}>> \\
(M_+\o (M_+\o M_-))\o M_- @>(1\o\beta_{23})\o 1>> \\
 (M_+\o (M_-\o M_+))\o M_- @>\text{associativity morphism}>> \\
(M_+\o M_-)\o (M_+\o M_-) @>v\o w>>V\o W,
\tag 2.1\endgather
$$
where $\beta_{23}$ denotes the braiding $\beta$ acting in the second 
and third components of the tensor product.

It is clear from this definition that all combinatorial complexity of the 
morphism $J$ comes from the arrows ``associativity morphism'' which involve 
associators. 

The arrows ``associativity morphism''
make the problem of checking various identities 
for $J$ (for example, the associativity identity) rather tedious.
To avoid this, we can use MacLane's theorem, which says that any monoidal 
category is equivalent to a strict one.  Namely,
when we check identities between morphisms in the category,
 we will assume that the category $\Cal M$ is replaced with 
an equivalent strict
monoidal category
and ignore associativity morphisms. 
For example, the definition of $J$ will look as follows:
$$
J_{VW}(v\o w)=(v\o w)\circ (1\o\beta_{23}\o 1)\circ(i_+\o i_-).
$$

However, when we do computations with vectors in modules from $\Cal M$,
it is important to pay attention to brackets, since different positions
of brackets are related with each other by the associator. 

\proclaim{Proposition 2.2} The maps $J_{VW}$ are isomorphisms
and define a tensor structure 
on the functor $F$.
\endproclaim

\demo{Proof.} 
It is obvious that $J_{VW}$ is an isomorphism since
it is an isomorphism modulo $h$. It is also clear that   
$J_{V\bold 1}=J_{\bold 1V}=1$. Thus
the only thing we need to check is the associativity identity
$J_{V\o W,U}\circ (J_{VW}\o 1)=J_{V,W\o U}\circ (1\o J_{WU})$.
To prove this equality, we need the following result.

\vskip .03in
{\bf Lemma 2.3} $(i_{\pm}\o 1)\circ i_{\pm}=(1\o i_{\pm})\circ i_{\pm}$
in $\text{Hom}(M_{\pm},M_{\pm}^{\o 3})$.
\vskip .03in

{\it Proof.} We prove the identity for $i_+$. 
The identity for $i_-$ is proved in the same way.

We need to show that for any vector $x\in M_{+}$ 
$$
\Phi\cdot (i_{+}\o 1) i_{+}x=(1\o i_{+})i_{+}x.\tag 2.2
$$
Since comultiplication in $U(\g_{-})$ is coassociative,
i.e. $(i_{+}\o 1) i_{+}x=(1\o i_{+}) i_{+}x$, 
it is sufficient to show that the associator $\Phi$ is the 
identity on the image
of $(i_{+}\o 1) i_{+}$. Because $\Phi$ is $\g$-invarint, 
it is enough to show that $\Phi\cdot (i_{+}\o 1) i_{+}1_+=
(i_{+}\o 1) i_{+}1_+$, i.e. 
$$
\Phi\cdot (1_+\o 1_+\o 1_+)=1_+\o 1_+\o 1_+.\tag 2.3
$$
Since the subalgebras $\g_+,\g_-$ are isotropic,  
the operators $\Omega_{12}$, $\Omega_{23}$ annihilate 
the vector $1_{+}\o 1_+\o 1_+$. 
Thus, equation (2.3) follows from the definition of $\Phi$.$\square$
\vskip .03in

Now we can finish the proof of the proposition. Let 
$\psi_1,\psi_2: M_+\o M_-\to (M_+\o M_-)^{\o 3}$ be the morphisms defined 
by
$$
\gather
\psi_1=(1\o \beta_{23}\o 1\o 1\o 1)\circ (i_+\o i_-\o 1\o 1)\circ 
(1\o\beta_{23}\o 1)\circ (i_+\o i_-),\\ 
\psi_2=(1\o 1\o 1\o \beta_{45}\o 1)\circ (1\o 1\o i_+\o i_-)\circ
(1\o\beta_{23}\o 1)\circ (i_+\o i_-),\tag 2.4
\endgather
$$
Then for any $v\in F(V),w\in F(W),u\in F(U)$
we have 
$$
\gather
J_{V\o W,U}(J_{VW}\o 1)(v\o w\o u)=(v\o w\o u)\circ \psi_1,\\
J_{V,W\o U}(1\o J_{WU})(v\o w\o u)=(v\o w\o u)\circ \psi_2.
\endgather
$$
Therefore, 
to prove the proposition, it is sufficient to show that $\psi_1=\psi_2$.

To prove this equality, we observe that the functoriality of the braiding
implies the identities
$$
\gather
(i_+\o i_-\o 1\o 1)\circ (1\o\beta_{23}\o 1)=(1\o\beta_{3,45}\o 1)
\circ (i_+\o 1\o i_-\o 1),\\
(1\o 1\o i_+\o i_-)\circ (1\o\beta_{23}\o 1)=(1\o\beta_{23,4}\o 1)
\circ (1\o i_+\o 1\o i_-)\tag 2.5\endgather
$$
(here $\beta_{3,45}$ means the braiding applied to the third factor and 
to the product of the fourth and the fifth factors). 
Using Lemma 2.3 and identities 
(2.5), we reduce the 
statement $\psi_1=\psi_2$ to the identity
$$
(1\o \beta_{23}\o 1\o 1\o 1)\circ (1\o\beta_{3,45}\o 1)=
(1\o 1\o 1\o \beta_{45}\o 1)\circ (1\o\beta_{23,4}\o 1),\tag 2.6
$$
which follows directly from the braiding axioms.
$\square$\enddemo

We call the functor $F$ equipped with the tensor structure $J$ the
fiber functor.

\vskip .1in
\centerline{\bf 3. Quantization of the double of a Lie bialgebra.}
\vskip .1in

3.1.{\it Topological Hopf algebras.}

Let $A$ be
an algebra over $k[[h]]$ with unit.
Let $I$ be a proper two-sided ideal in $A$ such that
$h\in I$. This ideal gives rise to a
translation invariant topology on $A$ such that
$\{I^n, n\ge 0\}$ is a basis of neighborhoods of $0$.
We will call $A$ a topological algebra if $A$
is complete in this topology, and $A/h^NA$ is a free 
$k[h]/(h^N)$-module for each $N\ge 1$. 

Let $A,B$ be two topological algebras, $I,J$ be the  
the corresponding ideals. Define $A\o B$ to be the projective limit of
algebras $A/I^n\o_{k[h]/h^n} B/J^n$ as $n\to\infty$.
Then $A\o B$ is also a topological algebra, with topology defined by
the ideal $I\o B+A\o J$.
 
We say that a topological algebra $A$ is a topological Hopf algebra if
it is equipped with comultiplication $\Delta: A\to A\o A$, 
the counit 
$\e:A\to k[[h]]$, 
and the
antipode $S:A\to A$, which are $k[[h]]$-linear, continuous, 
and satisfy the standard axioms of a Hopf algebra. 
Note that an infinite dimensional topological Hopf algebra may not be 
literally 
a Hopf algebra because the image of comultiplication may not belong to the 
algebraic tensor square of $A$.

Topological algebras and Hopf algebras over $k$ are defined similarly.

If $A$ is a topological algebra or Hopf algebra 
over $k[[h]]$ then $B=A/hA$ is
a topological algebra, respectively Hopf algebra, over $k$. 
In such a case we say that $A$ is a formal 
deformation of $B$. In particular, if $B=U(\g)$ with discrete topology, 
where $\g$ is a Lie 
algebra, then $A$ is called a quantized universal enveloping algebra
\cite{Dr1}. 

The following definition is due to Drinfeld \cite{Dr1}.

\proclaim{Definition} Let $(\g,\delta)$ be a Lie bialgebra. 
We say that a quantized universal enveloping algebra $A$ is a quantization of 
$(\g,\delta)$, or that $(\g,\delta)$ is the quasiclassical limit of $A$, if 

(i) $A/hA$ is isomorphic to $U(\g)$ as a Hopf algebra, and 

(ii) For any $x_0\in\g$ and any $x\in A$ equal to $x_0$ mod h one has
$$
h^{-1}(\Delta(x)-\Delta^{op}(x))\equiv \delta(x_0)\text{ mod h},
$$
where $\Delta^{op}$
is the opposite comultiplication ($\Delta^{op}=s\Delta$). 
\endproclaim
 
\vskip .05in
3.2. {\it The algebra of endomorphisms of
the fiber functor.}

Let $H=\text{End}(F)$ be the algebra of endomorphisms of the functor $F$.
This algebra is naturally isomorphic to $U(\g)[[h]]$. Namely, the map
$\alpha: U(\g)[[h]]\to H$ is defined on $x\in U(\g)$ by the formula
$(\alpha(x)f)(y)=f(yx)$, where $f\in \text{Hom}(U(\g),M)$, and 
is extended by linearity and continuity to $U(\g)[[h]]$.
This map is an isomorphism of algebras. From now on we will make no
distinction between $U(\g)[[h]]$ and $H$, identifying them by $\alpha$.

Let $F^2:\Cal M\times\Cal M\to \Cal A$ be the bifunctor defined by 
$F^2(V,W)=F(V)\o F(W)$. It is clear that $\text{End}(F^2)=H\o H$.

The algebra $H$ has a natural comultiplication $\Delta:H\to H\o H$
defined by $\Delta(a)_{V,W}(v\o w)=J_{VW}^{-1}a_{V\o W}
J_{VW}(v\o w)$, $a\in H$, $v\in F(V),w\in F(W)$
where $a_V$ denotes the action of $a$ in $F(V)$. We can also define the counit 
on $H$ by $\e(a)=a_{\bold 1}\in k[[h]]$,
where $\bold 1$ is the neutral object.

For any $V\in \Cal M$, let $V^*$ be the dual space to $V$
(regarded as an object of $\Cal M$), and let $\sigma_V:V^*\o V\to \bold 1$ 
be the canonical pairing. 
We have a functorial isomorphism 
$\xi_V: F(V^*)\to F(V)^*$
defined by $\xi_V(v^*)(v)=F(\sigma_V)J_{V^*V}(v^*\o v)$, $v\in F(V), 
v^*\in F(V^*)$. For any $a\in H$, let $\widetilde{S(a)}_V=
(\xi_V^*)^{-1}a_{V^*}^*\xi_V^*$ be a morphism $F(V)^{**}\to F(V)^{**}$.
It is easy to show that the subspace 
$F(V)\subset F(V)^{**}$ is invariant under 
this morphism. The 
antipode $S: H\to H$ is defined
by $S(a)_V=\widetilde{S(a)}_V|_{F(V)}$.   
  
\proclaim{Proposition 3.1}
The algebra $H$ equipped with $\Delta,\e,S$ is a 
topological Hopf algebra.
\endproclaim

The proof is straightforward.

\vskip .05in
 
3.3. {\it Explicit representation of comultiplication and antipode.}

Let $\Delta_0:U(\g)\to U(\g)\o U(\g)$ be the standard coproduct.
For any $V,W\in\Cal M$, let $J^0_{VW}: F(V)\o F(W)\to F(V\o W)$
be the morphism defined by the formula 
$J^0_{VW}(v\o w)(x)=(v\o w)(\Delta_0(x))$,
$x\in U(\g)$, $v\in F(V)$, $w\in F(W)$. It is clear
that $J_{VW}\equiv J^0_{VW}\text{ mod }h$.

Let $J\in U(\g)^{\o 2}[[h]]$ be defined by the formula
$$
J=
(\phi^{-1}\o\phi^{-1})\biggl(\Phi_{1,2,34}^{-1}(1\o \Phi_{2,3,4})
s e^{h\Omega_{23}/2}(1\o \Phi_{2,3,4}^{-1})\Phi_{1,2,34}
(1_+\o 1_+\o 1_-\o 1_-)\biggr),\tag 3.1
$$
where $\phi$ is the isomorphism 
of Lemma 2.1.

\proclaim{Proposition 3.2} For any $V,W\in\Cal A$, $v\in F(V)$, $w\in F(W)$
one has $J_{VW}(v\o w)=J_{VW}^0J(v\o w)$.
\endproclaim 

\demo{Proof} The statement follows from the definition (2.1) of $J_{VW}$.
$\square$\enddemo

\proclaim{Lemma 3.3} Let $a\in H$. Then
$$
\Delta(a)=J^{-1}\Delta_0(a)J.\tag 3.2
$$
\endproclaim

\demo{Proof} The lemma follows from Proposition 3.2 and the identities
$\Delta_0(a)_{V,W}=(J_{VW}^0)^{-1}a_{V\o W}J_{VW}^0$,
$\Delta(a)_{V,W}=J_{VW}^{-1}a_{V\o W}J_{VW}$ , $a\in U(\g)$.
$\square$\enddemo

Now consider the explicit expression for the antipode.
For any $V\in \Cal M$ define the morphism 
$\xi_V^0:F(V^*)\to F(V)^*$ by 
$\xi_V^0(v^*)(v)=F(\sigma_V)J_{V^*V}^0(v^*\o v)$,$v\in V$, $v^*\in V^*$. 
It is clear that $\xi_V\equiv \xi_V^0\text{ mod }h$.

Let $S_0:U(\g)\to U(\g)$ be the usual antipode.
Let $J=\sum_j x_j\o y_j$, $x_j,y_j\in U(\g)[[h]]$ 
(the sum is finite modulo $h^n$ for any $n$).
Define an element $Q\in U(\g)[[h]]$ by $Q=\sum_jS_0(x_j)y_j$.

\proclaim{Lemma 3.4} Let $a\in H$. Then 
$$
S(a)=Q^{-1}S_0(a)Q.\tag 3.3
$$
\endproclaim

\demo{Proof} It follows from the definitions of $\xi_V$, $\xi_V^0$,
and $Q$ that  $\xi_V=\xi_V^0S_0(Q)_{V^*}$. Thus the Lemma follows from the 
formulas $S(a)_V=(\xi_V^*)^{-1}a_{V^*}^*\xi_V^*|_{F(V)}$,
$S_0(a)_V=(\xi_V^{0*})^{-1}a_{V^*}^*\xi_V^{0*}|_{F(V)}$.
$\square$\enddemo

Thus, we have proved the following result.

\proclaim{Corollary 3.5} 
Introduce a new comultiplication and antipode on the topological Hopf 
algebra $U(\g)[[h]]$
by
$$
\Delta(x)=J^{-1}\Delta_0(x)J,\ S(x)=Q^{-1}S_0(x)Q,\tag 3.4
$$
where $\Delta_0,S_0$ are the usual comultiplication and antipode.
Then $(U(\g)[[h]],\Delta,S)$ is a topological Hopf algebra isomorphic to $H$.
\endproclaim

We will denote the topological Hopf algebra $(U(\g)[[h]],\Delta,S)$
by $U_h(\g)$.

{\bf Remark. } It is easy to see that according to the terminology of
\cite{Dr2}, the element $J^{-1}$ is a twist that realizes an equivalence
between the quasi-Hopf algebra $(U(\g)[[h]],\Phi)$ and the Hopf
algebra $U_h(\g)$. 
  
\vskip .05in
3.4.{\it The quasiclassical limit of $U_h(\g)$}.

\proclaim{Proposition 3.6} 
The topological Hopf algebra $U_h(\g)$ is a quantization of the Lie bialgebra
$(\g,\delta_\g)$.  
\endproclaim

\demo{Proof} Take $a\in\g\subset U_h(\g)$. 
Let $\delta(a)\in U(\g)\o U(\g)$ be defined by
the formula $\delta(a)=h^{-1}(\Delta(a)-\Delta^{op}(a))\text{ mod h}$. 
To prove the proposition, we need to show that for any $a\in\g$
one has $\delta(a)=\delta_\g(a)$, where $\delta_\g(a)$ is defined in
Chapter 1.

It is easy to check 
the following identities:
$$
e^{h\Omega/2}\equiv 1+h\Omega/2\text{ mod }h^2, \Phi\equiv 1\text{ mod }h^2.
\tag 3.5
$$

Let $\{g_j^+\}$ be a basis of $\g_+$, $\{g_j^-\}$ be 
the dual basis of $\g_-$, and $r=\sum_jg_j^+\o g_j^-$. 
Identities (3.1) and (3.5) imply that
$$
J\equiv 1+hr/2\text{ mod }h^2.\tag 3.6
$$
Therefore, by Lemma 3.3, 
$$
\Delta(a)\equiv \Delta_0(a)+\frac{h}{2}[\Delta_0(a),r]\text{ mod }h^2.\tag 3.7
$$
Thus,
$$
\Delta(a)-\Delta^{op}(a)\equiv \frac{h}{2}[\Delta_0(a),r-sr]\text{ mod }h^2.
\tag 3.8
$$
Since $r+sr$ ($=\Omega$) is $\g$-invariant, we obtain
$$
\delta=dr=\delta_\g,\tag 3.9
$$
Q.E.D.
$\square$\enddemo

\vskip .05in
3.5. {\it The quaitriangular structure on $U_h(\g)$}.

Define the element 
$$
R=(J^{op})^{-1}
e^{h\Omega/2}J\in U_h(\g)^{\o 2},\tag 3.10
$$ 
where $J^{op}$ is obtained 
from $J$ by permuting components. We call this element the universal 
R-matrix of $U_h(\g)$. 

\proclaim{Proposition 3.7} $R$ defines 
a quasitriangular structure on $U_h(\g)$.
That is, $R$ is invertible and
$$
R\Delta=\Delta^{op}R,\tag 3.11
$$
$$
 (\Delta\o 1)(R)=R_{13}R_{23},(1\o \Delta)(R)=
R_{13}R_{12}.\tag 3.12
$$ 
Moreover, $R$ is a quantization of $r$, i.e.  
$$
R\equiv 1+hr \text{ mod } h^2.\tag 3.13
$$
\endproclaim

\demo{Proof} Identity (3.13) follows from (3.5),(3.6)
and the definition of $R$. This identity implies that $R$ 
is invertible.

One has
$$
\gather
R\Delta(a)=(J^{op})^{-1}e^{h\Omega/2}J\Delta(a)=
(J^{op})^{-1}e^{h\Omega/2}\Delta_0(a)J=
(J^{op})^{-1}\Delta_0(a)e^{h\Omega/2}J=\\
\Delta^{op}(a)(J^{op})^{-1}e^{h\Omega/2}J=\Delta^{op}(a)R,\tag 3.14
\endgather
$$
which proves (3.11). 

Now let us prove the first identity of (3.12). The second identity
is proved analogously. 

According to the definition of $R$, for any $V,W\in \Cal M$, $v\in F(V)$,
$w\in F(W)$, one has $R(v\o w)=sJ_{WV}^{-1}F(\beta_{VW})J_{VW}$. Thus, 
for any $U\in \Cal M$, $u\in F(U)$ one has
$$
\gather
(\Delta\o 1)(R)(v\o w\o u)=(J_{VW}^{-1}\o 1)R(J_{VW}\o 1)(v\o w\o u)=\\
s_{12,3}(1\o J_{VW}^{-1})J_{U,V\o W}^{-1}F(\beta_{V\o W,U})
J_{V\o W,U}(J_{VW}\o 1)(v\o w\o u),\tag 3.15\endgather
$$
where $s_{12,3}$ is the permutation of the first two components with the 
third one.
Using the braiding property $\beta_{V\o W,U}=(\beta_{VU}\o 1)\circ (1\o
\beta_{WU})$, the associativity of $J_{VW}$, and the obvious identities
$J_{U\o V,W}^{-1}F(\beta_{VU}\o 1)J_{V\o U,W}=F(\beta_{VU})\o 1$,
$J_{V,U\o W}^{-1}F(1\o \beta_{WU})J_{V,W\o U}=1\o F(\beta_{WU})$,
one finds that the right hand side of (3.15) equals to 
$R_{13}R_{23}(v\o w\o u)$, as desired.
$\square$\enddemo

\vskip .1in
\centerline{\bf 4. Quantization of finite-dimensional Lie bialgebras.}
\vskip .1in

Our purpose in this section is to represent the quasitriangular topological Hopf algebra 
$U_h(\g)$ as a quantum double of another topological Hopf algebra, $U_h(\g_+)$.
The topological Hopf algebra $U_h(\g_+)$ will be a quantization of the Lie 
bialgebra $\g_+$.

\vskip .05in
 
4.1.{\it The algebras $U_h(\g_\pm)$.} 

As we have seen, the fiber functor $F$ which we used to construct the 
quantum group $U_h(\g)$, is represented by the object $M_+\o M_-$ of $\Cal M$.
Therefore, we have a homomorphism $\theta: \End(M_+\o M_-)\to \End(F)=U_h(\g)$ 
defined by $\theta(a)v=v\circ a$, $v\in F(V)$, $V\in \Cal M$, 
$a\in\End(M_+\o M_-)$.

\proclaim{Lemma 4.1} The map $\theta$ is an isomorphism.
\endproclaim

\demo{Proof} The Lemma follows from Lemma 2.1.
$\square$\enddemo

 Thus, we can identify $\End(M_+\o M_-)$ with $U_h(\g)$. From now on we make 
no distinction between them. 

Now let us define the subalgebras $U_h(\g_\pm)\subset U_h(\g)$. 

Let $x\in F(M_+)$.
Define the endomorphism $m_-(x)$ of $M_+\o M_-$ to be the composition
of the following morphisms in $\Cal M$: 
$m_-(x)=(x\o 1)\circ(1\o i_-)$. This defines a linear map
$m_-: F(M_+)\to U_h(\g)$. Denote the image of this map by $U_h(\g_-)$. 

Let $m_-^0(x)\in U(\g_-)$ be defined by the equation $x(1_+\o 1_-)=
m_-^0(x)1_+$.
It is easy to show that $m_-(x)\equiv m_-^0(x)\text{ mod }h$,
which implies that $m_-$ is an embedding.

A similar definition can be made for $x\in F(M_-)$. 
Define the endomorphism $m_+(x)$ of $M_+\o M_-$ to be the composition
of the following morphisms in $\Cal M$: 
$m_+(x)=(1\o x)\circ(i_+\o 1)$. This defines an injective linear map
$m_+: F(M_-)\to U_h(\g)$. 
Denote the image of this map by
$U_h(\g_+)$.

\proclaim{Proposition 4.2} $U_h(\g_\pm)$ are subalgebras in $U_h(\g)$.
\endproclaim

\demo{Proof} Let us give a proof for $U_h(\g_-)$. The proof for $U_h(\g_+)$ is 
analogous.

 Using Lemma 2.3, we obtain
$$
\gather
m_-(x)\circ m_-(y)=(x\o 1)\circ(1\o i_-)\circ (y\o 1)\circ(1\o i_-)=\\
=(x\o 1)\circ (y\o 1\o 1)\circ(1\o 1\o i_-)\circ(1\o i_-)=\\ 
=(x\o 1)\circ (y\o 1\o 1)\circ(1\o i_-\o 1)\circ(1\o i_-)=
(z\o 1)\circ (1\o i_-),\tag 4.1
\endgather
$$
where $z=x\circ (y\o 1)\circ (1\o i_-)\in F(M_+)$. 

So by the definition we get $m_-(x)\circ m_-(y)=m_-(z)$.
$\square$\enddemo

Note that the algebra $U_h(\g_-)$ is a deformation of the algebra
$U(\g_-)$. Indeed, we can define a linear isomorphism $\mu:U(\g_-)[[h]]
\to U_h(\g_-)$ by $\mu(a)(1_+\o 1_-)=a1_+$. This isomorphism
has the property $\mu(ab)=\mu(a)\circ \mu(b)\text{ mod }h^2$, 
which follows from (3.5), but in general $\mu(ab)\ne \mu(a)\circ \mu(b)$.

The subalgebra $U_h(\g_-)$ has a unit since it is a
deformation of the algebra with unit $U(\g_-)$. In fact, one can show that 
the unit equals to $\mu(1)$, $1\in U(\g_-)$.

Similar statements apply to the algebra $U(\g_+)$.

\proclaim{Proposition 4.3} The map $U_h(\g_+)\o U_h(\g_-)\to U_h(\g)$
given by $a\o b\to ab$ is an isomorphism.
\endproclaim

\demo{Proof} The statement is true because it holds modulo $h$.
$\square$\enddemo

\vskip .05in

4.2. {\it Polarization of the R-matrix.}

Define the element $\tilde R\in U_h(\g_+)\o U_h(\g_-)$ by
the identity 
$$
\tilde R\circ \beta^{-1}\circ (i_+\o i_-)=\beta\tag 4.2
$$
in $\text{Hom}(M_+\o M_-,M_-\o M_+)$. It is obvious that such
an element is unique. It can be computed as follows. 

Let $\nu: M_\pm[[h]]\to U_h(\g_\mp)$ be the linear isomorphism defined by
the equation $\nu(x(1_+\o 1_-))=m_\mp(x)$ for any $x\in F(M_\pm)$. Let
$K\in U(\g)^{\o 2}[[h]]$ be given by  
$$
K=
(\phi^{-1}\o\phi^{-1})\biggl(\Phi_{1,2,34}^{-1}(1\o \Phi_{2,3,4})
s e^{-h\Omega_{23}/2}(1\o \Phi_{2,3,4}^{-1})\Phi_{1,2,34}
(1_+\o 1_+\o 1_-\o 1_-)\biggr).\tag 4.3
$$
Then it is easy to check, using (4.2), that 
$$
\tilde R=(\nu\o \nu)(K^{-1}e^{h\Omega/2}(1_-\o 1_+)).\tag 4.4
$$

\proclaim{Proposition 4.4} $\tilde R=R$.
\endproclaim

\demo{Proof} According to (3.10), the R-matrix $R\in U_h(\g)\o U_h(\g)$ 
is defined by the condition that for any $V,W\in \Cal M$ and 
$v\in F(V),w\in F(W)$ one has the equality
$$
\gather
R^{op}(v\o w)\circ \beta_{23}\circ (i_+\o i_-)=\\
\beta\circ (w\o v)\circ \beta_{23}\circ (i_+\o i_-)\tag 4.5 
\endgather
$$
in $\text{Hom}(M_+\o M_-,V\o W)$.

By the functoriality of the braiding, $R^{op}(v\o w)=
\beta\circ R(w\o v)\circ \beta_{12,34}^{-1}$. Besides, 
$\beta_{12,34}=\beta_{23}\circ\beta_{12}\circ\beta_{34}\circ \beta_{23}$.
Substituting this into (4.5) and taking into account that 
$\beta\circ i_\pm=i_\pm$, we get
$$
\gather
R(w\o v)\circ \beta_{23}^{-1}\circ (i_+\o i_-)=\\
 (w\o v)\circ \beta_{23}\circ (i_+\o i_-)\tag 4.6 
\endgather
$$
in $\text{Hom}(M_+\o M_-,W\o V)$.

To show that $R=\tilde R$ we have to prove the identity
$$
\gather
(1\o \tilde R\o 1)
\circ (i_+\o 1\o 1\o i_-)\circ \beta_{23}^{-1}\circ (i_+\o i_-)=\\
\beta_{23}\circ (i_+\o i_-)\tag 4.7
\endgather
$$
in $\text{Hom}(M_+\o M_-,M_+\o M_-\o M_+\o M_-)$.

Interchanging the order of factors on the left hand side of (4.7) and
using Lemma 2.3, we can rewrite (4.7) in the form:
$$
\gather
(1\o \tilde R\o 1)
\circ \beta_{34}^{-1}\circ (1\o i_+\o i_-\o 1)\circ (i_+\o i_-)=\\
\beta_{23}\circ (i_+\o i_-)\tag 4.8 
\endgather
$$
in $\text{Hom}(M_+\o M_-,M_+\o M_-\o M_+\o M_-)$.

It is obvious that identity (4.8) follows from
is the definition of $\tilde R$.
The proposition is proved.
$\square$\enddemo

\vskip .05in
4.3. {\it Subalgebras $U_h(\g_\pm)$ in terms of the R-matrix.}

Let $U_h(\g_\pm)^*=\text{Hom}_{\Cal A}(U_h(\g_\pm),k[[h]])$.
Define $k[[h]]$-linear maps $\rho_\pm: U_h(\g_\mp)^*\to U_h(\g_\pm)$, by
$\rho_+(f)=(1\o f)(R)$, $\rho_-(f)=(f\o 1)(R)$. 
Let $U_\pm$ be the images of the maps $\rho_\pm$, and
$\tilde U_\pm$ be the closures of the $k[[h]]$-subalgebras generated by
$U_\pm$.

\proclaim{Proposition 4.5} 

$U_h(\g_\pm)
\o_{k[[h]]}k((h))$ is the h-adic completion of
$\tilde U_\pm\o_{k[[h]]}k((h))$.

\endproclaim

\demo{Proof} We prove the statement for $\tilde U_+$. 
The proof for $\tilde U_-$ is similar. 

We start with the following statement.

\vskip .03in
{\bf Lemma 4.6.}
For any $x\in U(\g_+)$ there exists an element $t_x\in \tilde U_+\o k((h))$ 
such that
$t_x=x+O(h)$. If $x$ has degree $\le m$ with respect to the standard 
filtration in $U(\g_+)$, then $t_x$ can be chosen in such a 
way that $h^mt_x\in \tilde U_+$.
\vskip .03in

{\it Proof of the Lemma.} It is clear that $1\in \tilde U_+$ 
since $1=\rho_+(\e)$. So we can set $t_1=1$. 

Now consider the case $x\in\g_+$.
Let $f\in U_h(\g_-)^*$ be any 
element such that $f(1)=0$ and $f(\tilde a)=\<x,a\>$ for any $a\in\g_-$
and $\tilde a\in U_h(\g_-)$ such that $\tilde a=a\text{ mod }h$.
Then it follows from (3.13) that $\rho_+(f)=hx+O(h^2)$. So we can let
$t_x=h^{-1}\rho_+(f)$. Thus, the Lemma is true for $x\in\g_+$.

Since $\tilde U_+$ is an algebra, the validity of the Lemma for
$x\in\g_+$ implies its validity for any $x\in U(\g_+)$. $\square$

Now we can prove the proposition.
Let $T_0\in U_h(\g_+)$. Let 
$x_0\in U(\g_+)$ be the reduction of $T_0$ 
$\text{ mod }h$. Then $T_0-t_{x_0}$ is divisible by $h$, 
so we can consider $T_1=h^{-1}(T_0-t_{x_0})$ and repeat our procedure.
This gives us a sequence $x_i\in U(\g_+)$, and $T_0=\sum_{m\ge 0}
t_{x_m}h^m $. This shows that $T_0$ belongs to the h-adic
completion of $\tilde U_+\o k((h))$, as desired.
$\square$\enddemo

\proclaim{Theorem 4.7} The subalgebras $U_h(\g_\pm)$ are Hopf 
subalgebras in $U_h(\g)$. 
\endproclaim

\demo{Proof} The fact that $U_h(\g_\pm)$ are closed under 
the comultiplication
$\Delta$ follows from Proposition 4.5 and identities (3.12).
The fact that $U_h(\g_\pm)$ are closed under the antipode $S$ follows from 
Proposition 4.5 and the identity $(S\o 1)(R)=R^{-1}$, which holds in any
quasitriangular Hopf algebra.
$\square$\enddemo

{\bf Remark.} In fact, it is possible to prove the following explicit formula 
for coproduct in $U_h(\g_\mp)$: for any $x\in F(M_\pm)$
$$
\Delta(m_\mp(x))=(m_\mp\o m_\mp)(J_{M_\pm M_\pm}^{-1}(i_\pm\circ x)).\tag 4.9
$$
The proof is a direct verification. 
A similar formula is contained in Proposition 9.3.

It is obvious that $U_h(\g_+)/hU_h(\g_+)$ is isomorphic to $U(\g_+)$
as a Hopf algebra. Therefore, $U_h(\g_+)$
is a quantized universal enveloping algebra.
It follows from Proposition 3.6 that its quasiclassical limit is the Lie
bialgebra $\g_+$. Similar statements apply to $U_h(\g_-)$.  

\vskip .05in
4.4. {\it Duality of quantized universal enveloping algebras and the quantum 
double.}

The following general constructions can be found in \cite{Dr1}.

If $A$ is a quantized universal enveloping algebra then the dual 
$A^*=\text{Hom}_{\Cal A}(A, k[[h]])$ carries a natural structure of
a topological algebra. Namely, for any $x,y\in A$, $f,g\in A^*$
$fg(x)=(f\o g)(\Delta(x))$, and the unit is $\e$. 
It can be shown that $A^*$ has a unique 
maximal ideal $I^*$, which is the kernel of the linear map $A\to  k$ given by
$f\to f(1)\text{ mod }h$. The topology on $A^*$ is defined by the condition 
that $\{(I^*)^n,n\ge 0\}$ is a basis of neighborhoods of zero.
This implies that the topological algebras $(A\o A)^*$ and $A^*\o A^*$ are 
isomorphic. 

The algebra $A^*$ has a natural structure of a topological 
Hopf algebra. Namely, the coproduct is defined by $\Delta(f)(x\o y)=f(xy)$, 
the counit is $1$, and the antipode is $S^*$.  
(The definition of coproduct makes sense since the algebra
$A^*\o A^*$ is isomorphic to 
$(A\o A)^*$).

As a topological 
$ k[[h]]$-module, $A^*$ is isomorphic to $ k[[X_1,...,X_N]][[h]]$.

Let $A$ be any quantized universal enveloping algebra.
Let $A^*$ be the dual algebra, and let $I^*$ be the maximal ideal
in $A^*$.
Consider the h-adic completion $A^\vee$
of the subalgebra $\sum_{n\ge 0} h^{-n}(I^*)^n$
in the algebra
$A^*\o_{k[[h]]} k((h))$. Then $A^\vee$ is a new quantized universal
enveloping algebra \cite{Dr1}. This algebra is called 
the dual quantized universal
enveloping algebra to $A$.

The algebra $A^*$ can be identified with a subalgebra in $A^\vee$ which is 
constructed as follows.

Let $\Delta^n:A\to A^{\o n}$ be the iterated coproduct maps:
$\Delta^0(a)=\e(a)$, $\Delta^1(a)=a$, $\Delta^2(a)=\Delta(a)$,
$\Delta^n(a)=(\Delta\o 1^{\o (n-2)})(\Delta^{n-1}(a))$, $n>2$. 

Let $\Sigma=\{i_1,...,i_k\}\subset \{1,...,n\}$, and $i_1<...<i_k$.
Let $j_\Sigma:A^{\o k}\to A^{\o n}$ be the homomorphism defined 
by $j_\Sigma(a_1\o...\o a_k)=b_1\o...\o b_n$, $a_1,...,a_k\in A$,
where $b_i=1$ if 
$i\notin\Sigma$, and $b_{i_m}=a_m$, $m=1,...,k$.
 
Let $\Delta_\Sigma(a)=j_\Sigma(\Delta^k(a))$, $a\in A$. 

Define linear mappings $\delta_n:A\to A^{\o n}$ for all $n\ge 1$
by
$$
\delta_n(a)=\sum_{\Sigma\subset \{1,...,n\}}(-1)^{n-|\Sigma|}\Delta_\Sigma(a)
$$
and a Hopf
subalgebra $A'=\{a\in A: \delta_n(a)\in h^nA^{\o n}\}$ in $A$. 

It is easy to check that $A^*=(A^\vee)'$.

If $A$ is any Hopf algebra, let $A^{op}$ denote 
the Hopf algebra $A$ with the comultiplication $\Delta$ replaced by 
$\Delta^{op}$, and the antipode $S$ replaced with $S^{-1}$. 
$A^{op}$ is also a Hopf algebra.

Now we can define the notion of the quantum double. 
Let $A$ be a quantized universal enveloping algebra.
Consider the $ k[[h]]$-module $D(A)=A\o (A^{\vee})^{op}$. 
Let $R\in A\o A^*\subset A\o (A^\vee)^{op}$ 
be the canonical element. We can regard $R$ as an element 
of $D(A)\o D(A)$ using the embedding $A\o (A^\vee)^{op}\to D(A)\o D(A)$
given by $x\o y\to x\o 1\o 1\o y$. 
Drinfeld \cite{Dr1} 
showed that there exists a unique structure of a topological 
Hopf algebra on $D(A)$ such that 

1) $A\o 1, 1\o (A^\vee)^{op}$ are Hopf subalgebras in $D(A)$,

2) $R$ defines a quasitriangular structure on $D(A)$, i.e. is invertible and 
satisfies (3.12), (3.13), and

3) The linear mapping $A\o (A^\vee)^{op}\to D(A)$ given by $a\o b\to ab$ 
is bijective. 

$D(A)$, equipped with this structure, is a quasitriangular quantized 
universal enveloping algebra. It is called the quantum double of $A$.
 
\vskip .05in
4.5. {\it The quantum double of $U_h(\g_+)$.}

\proclaim{Proposition 4.8} $\rho_+$ is a homomorphism of 
topological Hopf algebras $(U_h(\g_-)^{op})^*\to U_h(\g_+)$. 
$\rho_-$ is a homomorphism of 
topological Hopf algebras $U_h(\g_-)^*\to U_h(\g_+)^{op}$. 
\endproclaim

\demo{Proof} We only prove the first statement. The second one is proved 
analogously.

It is clear that $\rho_+$ is continuous. Also, for any
$f,g\in (U_h(\g_-)^{op})^*$ one has 
$$
\gather
\rho_+(fg)=(1\o fg)(R)=(1\o f\o g)((1\o\Delta^{op})(R))=(1\o f\o g)(R_{12}
R_{13})=\\
(1\o f)(R)\cdot (1\o g)(R)=\rho_+(f)\rho_+(g);\\
\Delta(\rho_+(f))=\Delta((1\o f)(R))=(1\o 1\o f)((\Delta\o 1)(R))=
\\   (1\o 1\o f)(R_{13}R_{23})=(1\o 1\o\Delta(f))(R_{13}R_{24})=
(\rho_+\o\rho_+)(\Delta(f)).\endgather
$$
It is obvious that $\rho_+(1)=1$ and $\e(\rho_+(f))=\e(f)$ for any $f$.
Also, it is easy to check that $\rho_+((S^{-1})^*f)=S(\rho_+(f))$.
The proposition is proved.
$\square$\enddemo

\proclaim{Corollary 4.9} $U_\pm$ are Hopf subalgebras in $U_h(\g_\pm)$.
In particular, $\tilde U_\pm=U_\pm$.
\endproclaim

\demo{Proof} The first statement is clear. The second statement follows 
from the first one and the fact that $U_\pm$ is closed in $U_h(\g_\pm)$, 
which is easy to check.
$\square$\enddemo

\proclaim{Proposition 4.10} The maps $\rho_+$, $\rho_-$ are injective. 
\endproclaim

\demo{Proof} We show the injectivity of $\rho_+$ (the case of
$\rho_-$ is similar). Fix an element
 $f\in U_h(\g_-)^*,f\ne 0$. We can always assume that
$f\ne 0\text{ mod }h$. Let $x\in U(\g_-)$ be such that 
$f(t_x)\ne 0\text{ mod h}$ (where $t_x$ was defined in Lemma 4.6), 
$n\ge 0$ be such that $h^nt_x\in U_-$, and 
$g\in U_h(\g_+)^*$ be such that $\rho_-(g)=h^nt_x$. Such a $g$ exists 
by the definition of $n$. Then $g(\rho_+(f))=(g\o f)(R)=f(\rho_-(g))=
h^nf(t_x)\ne 0$.
Therefore, $\rho_+(f)\ne 0$.
$\square$\enddemo

\proclaim{Proposition 4.11} 
$U_\pm=U_h(\g_\pm)'$.
\endproclaim
  
\demo{Proof} We give the proof for $U_+$. The proof for $U_-$ is similar.
 
First we need the following statement.

\vskip .03in
{\bf Lemma 4.12} Let $t\in U_h(\g_+)'$ be an element such that
 $h^{-n}t\in U_h(\g_+)$ and
$h^{-n}t=x+O(h)$, $x\in U(\g_+)$, $x\ne 0$. 
Then $x$ has degree $\le n$. 
\vskip .03in

{\it Proof of the Lemma.} By the definition, $\delta_{n+1}(h^{-n}t)$ is 
divisible by $h$. On the other hand, $\delta_{n+1}(h^{-n}t)=
\delta_{n+1}(x)+O(h)$. Thus, $\delta_{n+1}(x)=0$, which implies that 
the degree of $x$ is $\le n$, since the kernel of $\delta_{n+1}$ on $U(\g_+)$
is the set of all
elements of $U(\g_+)$ whose degree is $\le n$. $\square$ 

Now we can prove the proposition.
By Lemma 4.6, for any $x\in U(\g_+)$ 
of degree $\le n$, an element $t_x$ can be chosen in such a way that 
$h^nt_x\in U_+$. This implies the inclusion $U_+\supset U_h(\g_+)'$.
Indeed, let
$T_0\in U_h(\g_+)'$, and $T_0\equiv h^mx_0\text{ mod }h^{m+1}$, where
$x_0\in U(\g_+)$. Then, according to Lemma 4.12, 
the degree of $x_0$ is $\le m$. Therefore, $h^mt_{x_0}\in U_+$. Thus,
 $T_1=T_0-h^mt_{x_0}\in U_+$ and is divisible by $h^{m+1}$, 
so we can repeat our procedure.
This gives us a sequence of elements 
$x_i\in U(\g_+)$ of degrees $m_i$ ($m_0=m$), such that $m_0<m_1<...<m_i<...$,
and $T_0=\sum_{i\ge 0}
t_{x_i}h^{m_i}$. This shows that $T_0$ belongs to $U_+$,
as desired.

To demonstrate the inclusion $U_+\subset U_h(\g_+)'$, 
observe that according to (3.12), 
$$
(\Delta^n\o 1)(R)=R_{1n+1}...R_{nn+1}.
$$
This implies that 
$$
(\delta_n\o 1)(R)=(R_{1n+1}-1)...(R_{nn+1}-1)=O(h^n).
$$
Therefore, $\delta_n(\rho_+(f))$ is divisible by $h^n$ for any
$f\in U_h(\g_-)^*$. 
$\square$\enddemo 

Comparing our results with the definitions of the previous, section
we see that we have obtained the following result.

\proclaim{Theorem 4.13} Let $\g_+$ be a finite dimensional Lie bialgebra
and $(\g,\g_+,\g_-)$ be the associated Manin triple. Then 

(i) There exist
quantized universal enveloping algebras $U_h(\g)$ and 
$U_h(\g_\pm)\subset U_h(\g)$, 
which are quantizations of the Lie bialgebras $\g$, $\g_\pm\subset \g$,
respectively; 

(ii) The multiplication map $U_h(\g_+)\o U_h(\g_-)\to U_h(\g)$ is a 
linear isomorphism;

(iii) 
The algebras $U_h(\g_+),U_h(\g_-)^{op}$ are dual to each other as 
quantized universal enveloping algebras, in the sense of Drinfeld \cite{Dr1};

(iv) The factorization $U_h(\g)=U_h(\g_+)U_h(\g_-)$ defines an isomorphism
of $U_h(\g)$ with the quantum double of $U_h(\g_+)$;

(v) $U_h(\g)$ is isomorphic to $U(\g)[[h]]$ as a topological algebra.
\endproclaim

\vskip .05in
\centerline{\bf 5. Quantization of solutions of the classical
Yang-Baxter equation.}
\vskip .05in

Let $A$ be an associative algebra over $ k$ with unit,
and $r\in A\o A$. The element $r$ is called a classical r-matrix if it satisies
the classical Yang-Baxter equation
$$
[r_{12},r_{13}]+[r_{12},r_{23}]+[r_{13},r_{23}]=0.\tag 5.1
$$
We say that $r$ is unitary if $r^{op}=-r$. An algebra $A$ equipped
with a classical $r$-matrix $r$ is called a classical Yang-Baxter
algebra. $A$ is called unitary if $r$ is unitary.

Let $A$ be a topological algebra over $k[[h]]$.
Let $R\in A\o A$. We say that $R$ is a quantum $R$-matrix if 
it satisfies the quantum Yang-Baxter equation
$$
R_{12}R_{13}R_{23}=R_{23}R_{13}R_{12}\tag 5.2
$$
We say that $R$ is unitary if $R^{op}=R^{-1}$. 
A topological algebra $A$ equipped
with a quantum $R$-matrix $R$ is called a quantum Yang-Baxter
algebra. $A$ is called unitary if $R$ is unitary.

The following theorem answers question 3.1 in \cite{Dr3}.
It shows that any classical Yang-Baxter algebra can be quantized. 

\proclaim{Theorem 5.1} Let $A$ be an associative algebra with 
unit over $k$, and
$r\in A\o A$ be a classical $r$-matrix. Then there exists a quantum $R$-matrix
$R\in A\o A[[h]]$ such that $R=1+hr\text{ mod }h^2$.
If in addition $r$ is unitary then $R$ can also be chosen
unitary.
\endproclaim

{\it Proof.}

We start with a construction of Reshetikhin and Semenov-Tian-Shansky
\cite{RS}.
Let $\g_+=\{(1\o f)(r),f\in A^*\},\g_-=\{(f\o 1)(r),f\in A^*\}$
be vector subspaces in $A$. It is clear that $\g_+,\g_-$ are 
finite-dimensional, $r\in\g_+\o \g_-$, and 
the map $\chi_r:\g_+^*\to\g_-$ defined by $\chi_r(f)=(f\o 1)(r)$, is an
isomorphism of vector spaces.

{\bf Remark.}
Note that the spaces $\g_+$ and $\g_-$ may intersect nontrivially and even 
coincide. 

\proclaim{Lemma 5.2.} $\g_+,\g_-$ are Lie subalgebras in $A$. 
\endproclaim

\demo{Proof} Let $x,y\in\g_+$, $x=(1\o f)(r), y=(1\o g)(r)$. Using (5.1), we have
$$
[xy]=(1\o f\o g)([r_{12}r_{13}])=-(1\o f\o g)([r_{12}+r_{13},r_{23}])=
(1\o h)(r),\tag 5.3
$$
where $h\in A^*$, $h(a)=(f\o g)([r,a\o 1+1\o a])$. Thus, $[xy]\in\g_+$,
i.e. $\g_+$ is a Lie algebra. The proof for $\g_-$ is similar.
$\square$\enddemo

Let $\g=\g_+\oplus\g_-$ be a vector space 
Define the skew-symmetric bracket $[,]:\g\o\g\to\g$ as follows.
If $x,y\in\g_+$ or $x,y\in\g_-$ then the bracket $[xy]$ is the Lie bracket
in $\g_+$ or $\g_-$, respectively. If $x\in\g_+,y\in\g_-$, then
$[xy]$ is defined by 
$$
[xy]=(\ad^*x)y-(ad^*y)x. \tag 5.4
$$

Let $\pi:\g\to A$ be the linear map whose restrictions to 
$\g_+,\g_-$ are the corresponding embeddings.
The restrictions of $\pi$ to $\g_+$, $\g_-$ are injective
but in general $\pi$ itself is not an embedding.

\proclaim{Lemma 5.3} $\pi([xy])=[\pi(x),\pi(y)]$, $x,y\in\g$.
\endproclaim

\demo{Proof} The Lemma is a direct consequence of the classical 
Yang-Baxter equation.
$\square$\enddemo

\proclaim{Lemma 5.4} $(\g,[,])$ is a Lie algebra.
\endproclaim

\demo{Proof} We have to check the Jacobi identity in $\g$. 
It is enough to check it for three elements $a,x,y$ such that 
$a\in\g_+$, $x,y\in\g_-$. For brevity we write $a(x)$ for $(\ad^*a)x$.
We have
$$
\gather
[a[xy]]=a([xy])-[xy](a),\\
[y[ax]]=[y,a(x)-x(a)]=[y,a(x)]-y(x(a))+y(a(x)),\\
[x[ya]]=[x,y(a)-a(y)]=-[x,a(y)]+x(y(a))-x(a(y)).\tag 5.5
\endgather
$$
Adding these three identities, and using the fact that 
$[xy](a)=x(y(a))-y(x(a))$, we get
$$
[a[xy]]+[y[ax]]+[x[ya]]=a([xy])+[y,a(x)]-[x,a(y)]+y(a(x))-x(a(y)).
\tag 5.6
$$

Denote the right hand side of (5.6) by $X$.
Applying $\pi$ to both sides of (5.6), and using Lemma 5.3
 and the Jacobi identity in $A$,
we get
$$
\pi(X)=0.
\tag 5.7
$$
Since $X\in\g_+$, and $\pi$ is injective on $\g_+$, we get $X=0$,
which implies the Jacobi identity in $\g$.
$\square$\enddemo

Let $\<,\>$ be the inner product on $\g$
such that $\<x_++x_-,y_++y_-\>=x_-\cdot y_++y_-\cdot x_+$, where 
$x_+,y_+\in\g_+$, $x_-,y_-\in\g_-$, and the dot denotes the natural pairing 
$\g_-\o\g_+\to k$ defined by the map $\chi_r$. 
This inner product is ad-invariant.
Thus, $(\g,\g_+,\g_-)$ is a Manin triple.

Now we can finish the proof of the theorem. Lemma 5.3 implies that 
$\pi:\g\to A$ is a homomorphism of Lie algebras. Therefore, it extends to 
a homomorphism of associative algebras $\pi:U(\g)\to A$. 
Furthermore, $(\g,\g_+,\g_-)$ is a Manin triple.
The Lie bialgebra $\g$ is quasitriangular, and its quasitriangular structure
is defined by the classical r-matrix $\tilde r=\sum x_+^i\o x_-^i$, where
$x_+^i$ is a basis of $\g_+$, and $x_-^i$ is a dual basis of $\g_-$.
Note that $(\pi\o\pi)(\tilde r)=r$.
  
By Theorem 4.13, there exists a quasitriangular topological Hopf algebra 
$U_h(\g)$, with a quasitriangular structure $\tilde R\in U_h(\g)\o U_h(\g)$.
Moreover, the associative algebra $U_h(\g)$ is isomorphic to $U(\g)[[h]]$,
and the isomorphism can be chosen to be the identity modulo $h$. Thus, we can 
assume that $\tilde R\in (U(\g)\o U(\g))[[h]]$.

Set $R=(\pi \o\pi)(\tilde R)$. From what we said above it follows that
$R$ satisfies (5.2) and $R=1+hr$ modulo $h^2$.

Assume now that $r^{op}=-r$. Let $\tilde \Omega=\tilde r+\tilde r^{op}$.
It follows immediately
from the construction of $\tilde R$ that $\tilde R^{op}\tilde
R$ is conjugate to 
$e^{h\tilde \Omega}$. But $(\pi\o\pi)(\tilde\Omega)=r+r^{op}=0$.
This implies that $R^{op}R=1$, as desired.

The theorem is proved.$\square$

Let $\Cal R$ be 
the ring of algebraic functions of a variable $h$ with coefficients in $k$ 
which are regular at $h=0$. 

\proclaim{Theorem 5.5} Let $A$ be a finite-dimensional associative algebra 
with unit 
over $k$, and $r\in A\o A$ be a classical $r$-matrix. 
Then there exists a family of quantum R-matrices
$R(h)\in A\o A\o \Cal R$ such that $R=1+hr+O(h^2)$, $h\to 0$. 
If in addition $r$ is unitary then $R(h)$ can also be chosen
unitary.
\endproclaim

\demo{Proof} The theorem follows immediately from Theorem 5.1 
and the following result of M.Artin 
\cite{Ar}.

{\bf Theorem.} Any system of polynomial equations in indeterminates 
$x_1,...,x_n$ with coefficients in $k[h]$ which has solutions over $k[[h]]$
also has solutions over $\Cal R$.

Indeed, let us write $R$ in the form $R=1+hr+h^2X(h)$, and look for a series
$X(h)$ such that $R$ satisfies the quantum Yang-Baxter equation, and the 
unitarity condition in the case when $r$ is unitary. This is a system of 
polynomial equations on the components of $X(h)$ with coefficients in $k[h]$.
By Theorem 5.1, it has solutions over $k[[h]]$. Therefore, by Artin's theorem,
it has solutions over $\Cal R$.  
\enddemo

\vskip .1in
\centerline {\bf 6. Quantization of quasitriangular
Lie bialgebras.}
\vskip .1in

6.1. {\it Quasitriangular quantization of quasitriangular Lie bialgebras.}
In this section we give a recipe of quantization of a quasitriangular 
Lie bialgebra $\a$ (not necessarily finite dimensional), 
which produces a quantized universal enveloping algebra
isomorphic to $U(\a)[[h]]$ as a topological algebra. 
This answers questions from Section 4 of \cite{Dr3}.

Let $\g_+=\{(1\o f)(r),f\in\a^*\}$, $\g_-=\{(f\o 1)(r),f\in\a^*\}$.
be subspaces in $\a$. By Lemma 5.2, applied to $A=U(\a)$,
these subspaces are finite dimensional Lie subalgebras 
in $\a$. Moreover, let $\g$ be the vector space $\g_+\oplus\g_-$.
This space is a Lie algebra with bracket defined by (5.4) and an invariant
inner product. By Lemma 5.3, we have a natural homomorphism
of Lie algebras $\pi:\g\to \a$, and it is easy to see that this homomorphism 
is a morphism of quasitriangular Lie bialgebras. 

Let $\Cal M_\a$ be the category whose
objects are $\a$-modules, and morphisms are defined by 
$\Hom_{\Cal M_\a}(V,W)=\Hom_\a(V,W)[[h]]$. 
 Let $\Cal M_\g$
be the Drinfeld category associated to $\g$. We have the pullback functor 
$\pi^*:\Cal M_\a\to\Cal M_\g$. 
Define the braided monoidal structure on $\Cal M_\a$ 
to be the pullback of the braided monoidal structure on $\Cal M_\g$.
This definition makes sense, since the element
$\Omega=r+r^{op}\in\g\o\g$
is $\g$-invariant by the definition of a quasitriangular Lie
bialgebra.

Let $M_+,M_-$ be the Verma modules 
over $\g$. Define a functor $F: \Cal M_\a\to \Cal A$ by 
$F(V)=\Hom_{\Cal M_\g}(M_+\o M_-,\pi^*(V))$. 
The tensor structure on $F$ is introduced in the same way as in
Section 1.8.
Let $H=\End F$. Since the functor $F$ is isomorphic to
the ``forgetful'' functor $V\to ``\text{ the }k[[h]]\text{ module
}V[[h]]$'', the algebra
 $H$ is isomorphic to $U(\a)[[h]]$ as a 
topological algebra over $k[[h]]$. On the other hand, $H$ has a natural coproduct 
and antipode defined analogously to Section 3.2, and a quasitriangular structure 
$R\in H\o H$ defined analogously to Section 3.5. It is easy to check that   
the quasiclassical limit of $H$ is the Lie bialgebra $\a$, and $R=1+hr+O(h^2)$, so
$r$ is the quasiclassical limit of $R$. 

Furthermore, suppose that the original Lie bialgebra $r$ is triangular, i.e. $r$ is 
a unitary $r$-matrix. Then $\Omega=r+r^{op}=0$, and hence $R^{op}R=J^{-1}e^{h\Omega}J=1$,
so the Hopf algebra $H$ is triangular, too.

Thus, we have the following 
theorem.

\proclaim{Theorem 6.1} 
Any quasitriangular Lie bialgebra $\a$ admits a quantization $U_h(\a)$ which is
a quasitriangular quantized universal enveloping algebra isomorphic to $U(\a)[[h]]$
as a topological algebra. If $\a$ is triangular, so is $U_h(\a)$. 
\endproclaim

\vskip .05in
6.2. 
{\it Identification of two quantizations of a quasitriangular Lie bialgebra.}
Let $\a$ be a finite dimensional
quasitriangular Lie bialgebra. Let $U_h(\a)$ be the quantization 
of $\a$ constructed in Section 4, and $U_h^{qt}(\a)$ be the quasitriangular 
quantization of $\a$ constructed in Section 6.1.

\proclaim{Theorem 6.2} The quantized universal enveloping algebras $U_h(\a)$, 
$U_h^{qt}(\a)$ are isomorphic.
\endproclaim

The proof of this theorem is given below and uses the functoriality
of quantization, which is proved in Chapter 10.

\proclaim{Corollary 6.3} The quantization of the double $\g$ of a 
finite dimensional Lie bialgebra $\a$ constructed in Chapter 3 
is isomorphic to the quantization of $\g$ as a Lie bialgebra, constructed in 
Chapter 4.
\endproclaim
  
To prove Theorem 6.2, we first need the following result, which appears
(in somewhat different form) in \cite{RS}.

\proclaim{Lemma 6.4} Let $\a$ be a quasitriangular Lie bialgebra, and $\g$ be the double 
of $\a$. Then the linear map $\tau:\g\to\a$ defined by
$$
\tau(x+f)=x+(f\o 1)(r), x\in\a, f\in\a^*,\tag 6.1
$$
is a homomorphism of quasitriangular Lie bialgebras.
\endproclaim

\demo{Proof} First we show that $\tau$ is a homomorphism of Lie algebras,
i.e. $\tau([g_1g_2])=[\tau(g_1)\tau(g_2)]$. 
This is obvious when $g_1,g_2\in\a$.
Assume that $f,g\in\a^*$. Then, using the classical Yang-Baxter equation, we get
$$
\gather
\tau([fg])=([fg]\o 1)(r)=(f\o g\o 1)((\delta\o 1)(r))=\\
(f\o g\o 1)([r_{13}+r_{23},r_{12}])=(f\o g\o 1)([r_{13},r_{23}])=[\tau(f)\tau(g)].
\tag 6.2\endgather
$$  
Now assume that $x\in\a,f\in \a^*$. Then
$$
\gather
\tau([xf])=\tau(\ad^*x(f))-\tau(\ad^*f(x))=\\
\tau((f\o 1)([r,x\o 1]))+
\tau((f\o 1)([x\o 1+1\o x,r]))=\\
\tau((f\o 1)([1\o x,r]))=[\tau(x)\tau(f)].\tag 6.3
\endgather
$$

Now we check that $\tau$ is a homomorphism of 
quasitriangular Lie bialgebras. Let $\tilde r$ be the quasitriangular structure on $\g$.
If $x_i$ is a basis of $\a$, and $f_i$ is the dual basis of $\a^*$, then 
$\tilde r$ is given by the formula $\tilde r=\sum_i x_i\o f_i$. 
Thus we have
$$
(\tau\o\tau)(\tilde r)=\sum_i\tau(x_i)\o\tau(f_i)=\sum_i x_i\o (f_i\o 1)(r)=r.\tag 6.4
$$
The Lemma is proved.
$\square$\enddemo

\demo{Proof of Theorem 6.2} Lemma 6.4 claims that there exists a 
morphism of quasitriangular Lie bialgebras $\tau:\g\to\a$ which is 
the identity on $\a$. Theorem 10.6 below states that 
quasitriangular quantization of Section 6.1 is a functor from the category
of quasitriangular Lie bialgebras to the category of quasitriangular
topological Hopf algebras over $k[[h]]$. Thus, $\tau$ defines a morphism
$\hat\tau: U_h^{qt}(\g)\to U_h^{qt}(\a)$. On the other hand, $U_h(\a)$ was constructed
as a subalgebra in $U_h^{qt}(\g)$, so we have 
an embedding $\eta: U_h(\a)\to U_h^{qt}(\g)$. Consider the morphism
$\tau\circ \eta: U_h(\a)\to U_h^{qt}(\a)$. This morphism is an isomorphism since it equals
to $1$ modulo $h$. The theorem is proved.
$\square$\enddemo

{\bf Remark.} An analogous theorem holds for infinite dimensional Lie bialgebras. 
Namely, the ``usual'' quantization of $\a$ defined in Section 9 is isomorphic 
to its quasitriangular quantization. The proof is analogous to the finite 
dimensional case.  

\vskip .05in

6.3. {\it Representations of $U_h(\g)$.}

Let $\a$ be a quasitriangular Lie bialgebra (not necessarily
finite dimensional). 
By a representation of $U_h(\a)$ we mean a topologically free 
$k[[h]]$-module $V$ together with a homomorphism 
$\pi_V: U_h(\a)\to \text{End}_{k[[h]]}V$. Representations of
$U_h(\g)$ form a braided tensor category, with the trivial 
associativity morphism and braiding 
defined by the $R$-matrix. Denote this category by 
$\Cal R$. 

The functor $F:\Cal M_{\a}\to\Cal A$ can be regarded as a functor 
from $\Cal M_{\a}$ 
to $\Cal R$, since for any $W\in \Cal M_\a$ the $k[[h]]$-module 
$F(W)$ is equipped with a natural action of $U_h(\g)$. 
We denote this new functor also by $F$. 
This functor inherits the
 tensor structure
defined by the maps $J_{VW}$. 

\proclaim{Theorem 6.5} The functor $F$ defines an equivalence of braided tensor
categories $\Cal M_\a\to \Cal R$. 
\endproclaim

\demo{Proof} The theorem follows from the definition of
the functor $F$, the algebra $U_h(\g)$ and the R-matrix $R$. 
$\square$\enddemo

\vskip .1in
\centerline{\bf Part II}
\vskip .1in
\centerline{\bf 7. Drinfeld category for an arbitrary Lie bialgebra.} 
\vskip .1in

7.1. {\it Topological vector spaces.} 
Recall the definition of the product topology.
Let $S$ be a set, $T$ a topological space, and 
$T^S$ the space of functions from $S$ to $T$. This space has 
a natural weak topology, which is the weakest of the topologies in which 
all the evaluation maps $T^S\to T$, $f\to f(s)$, are continuous. 
Namely, let $B$ be a basis of the topology on $T$. For any integer 
$n\ge 1$, elements 
$s_1,...,s_n\in S$, and open sets $U_1,...,U_n\in B$, define 
$V(s_1,...,s_n,U_1,...,U_n)=\{f\in T^S: f(s_i)\in U_i,i=1,...,n\}$. 
Let $\Cal B$ be the collection of all such sets $V$. This is a basis of a 
topology on $T^S$ which is called the weak topology. The obtained 
topological space is the product of copies of $T$ corresponding to 
elements of $S$. If $X$ is any subset in $T^S$, the weak topology on 
$T^S$ induces a topology on $X$. We will call it the weak topology as well. 

Let $ k$ be a field of characteristic zero with the discrete topology.
Let $V$ be a topological vector space over $k$. The topology on $V$
is called linear if open subspaces of $V$ form a basis of neighborhoods of 
$0$. 

{\bf Remark.} It is clear that in any topological vector space, an
open subspace is also closed. 
 
Let $V$ be a topological vector space over $k$ with linear topology.
$V$ is called separated if the map $V\to \underleftarrow{\lim}(V/U)$ 
is a monomorphism, where $U$ runs over open subspaces of $V$. 

Topology on all vector spaces we consider in this paper 
will be linear and separated. 
so we will say ``topological vector space'' for 
``separated topological vector space 
with linear topology''.
    
Let $M$, $N$ be topological vector spaces over $ k$.
We denote by $\Hom_ k(M,N)$ 
the space of continuous linear operators from $M$ to $N$, 
equipped with the weak topology.  
It is clear that a basis of neighborhoods of zero in
$\Hom_k(M,N)$ is generated by sets of the form
$\{A\in \Hom_k(M,N): 
Av\in U\}$, where $v\in M$, and $U\subset N$ is an open set. 

In particular, if $N= k$ with the discrete topology,
the space $\Hom_ k(M,N)$ is the space of all continuous linear functionals 
on $M$, which we denote by $M^*$.
It is clear that a basis of neighborhoods of zero
in $M^*$ consists of
 orthogonal complements of finite-dimensional subspaces in $M$.
In particular, if $M$ is discrete then the canonical embedding
$M\to (M^*)^*$ is an isomorphism of linear spaces. However, if
$M$ is infinite-dimensional, this embedding is not an isomorphism
of topological vector spaces since the space $(M^*)^*$ is not discrete.
\vskip .05in

7.2 {\it Complete vector spaces}

Let $V$ be a topological vector space over $k$.
$V$ is called complete if the map $V\to \underleftarrow{\lim}(V/U)$ 
is a epimorphism, where $U$ runs over open subspaces of $V$. 

In particular, if a complete space $M$ has a countable
basis of neighborhoods of $0$, then there exists a 
filtration $M=M_0\supset M_1\supset ...$, such that $ \cap_{n\ge 0} M_n=0$,
and $\{M_n\}$ is a basis of neighborhoods of zero in $M$. 
In this case $M=\lim_{n\to \infty} M/M_n$. 

{\it Examples.} 1. Any discrete vector space is complete.

2. If $V$ is a discrete vector space then the 
topological space $M=V[[h]]$
of formal power series in $h$ with coefficients in $V$ is a complete
vector space.

Let $V$ be a complete vector space, $U\subset V$ an open subspace. 
Then $U$ is complete and $V/U$ is discrete. 

Let $V,W$ be complete vector spaces. Consider the space
$V\ho W=\underleftarrow{\lim}V/V_1\o W/W_1$, where the projective
limit is taken over open subspaces $V_1\subset V$, $W_1\subset W$. 
It is easy to see that $V\ho W$ is a complete vector space.
We call the operation $\ho$ the completed tensor product.

A basis of neighborhoods of $0$ in $V\ho W$ is the collection
of subspaces $V\ho W_1+V_1\ho W$, where $V_1,W_1$ are open subspaces 
in $V,W$.  

{\it Example.} Let $V$ be a discrete space. Then $V\ho k[[h]]=V[[h]]$.

Complete vector spaces form an additive category in which morphisms are 
continuous linear operators. This category, equipped with 
tensor product $\ho$, is a strict symmetric tensor 
category.

\vskip .05in
7.3. {\it Equicontinuous $\g$-modules.}

Let $M$ be a topological vector space over $k$, and
 $\{A_x,x\in X\}$ be a family of elements of $\End M$. 
We say that the family 
$\{A_x\}$ is equicontinuous if for every neighborhood
of the origin $U\subset M$ there exists another neighborhood of the origin
$U'\subset M$ such that $A_xU'=U$ for all $x\in X$. 
For example, if $M$ is complete and
$A\in \End M$ is any continuous linear operator,
then $\{\lambda A,\lambda\in k\}$ is equicontinuous. 

Fix a topological Lie algebra $\g$.

{\bf Definition.}
Let $M$ be a complete vector space.
We say that $M$ is an equicontinuous 
$\g$-module if one is given a continuous homomorphism
of topological Lie algebras $\pi:\g\to\End M$, such that the family of 
operators $\pi(g), g\in\g$, is equicontinuous. 

{\it Example.} If $M$ is a complete vector space 
with a trivial $\g$-module structure then $M$ is an equicontinuous 
$\g$-module.

Let $V,W$ be equicontinuous $\g$-modules. It is easy to check that
 $V\ho W$ has a natural structure of 
an equicontinuous $\g$-module. Moreover, $(V\ho W)\ho U$ is naturally identified with
$V\ho (W\ho U)$ for any equicontinuous $\g$-modules $V,W,U$. 
This means that the category of equicontinuous $\g$-modules, where morphisms
are continuous homomorphisms, is a monoidal category. This category 
is symmetric since 
the objects $V\ho W$ and $W\ho V$ are identified by the 
permutation of components. We denote this category by $\Cal M^e_0$.

\vskip .05in  
7.4. {\it Lie bialgebras and Manin triples.}
Let $\a$ be a Lie bialgebra over $ k$.
We will regard $\a$ as a topological Lie algebra with the discrete topology.
Let $\a^*$ be the full dual space to $\a$.
The cocommutator defines a Lie bracket on
$\a^*$ which is continuous in the weak topology, so $\a^*$
has a natural structure of a topological Lie algebra.   

Furthermore, the space $\a\oplus\a^*$ has a natural topology, and
the Lie bracket on $\g$ defined
by (1.1) is continuous in this topology.

Let $\g$ be a 
Lie algebra with a nondegenerate invariant
inner product $\<,\>$. So far we have no topology on $\g$. 
Let $\g_+$,$\g_-$ be isotropic Lie subalgebras in $\g$,
such that $\g=\g_+\oplus\g_-$ as a vector space.
The inner product $\<,\>$ defines an embedding $\g_-\to\g_+^*$.
If this embedding is an isomorphism then we equip $\g$ with
a topology, by putting the discrete topology on $\g_+$ and the weak 
topology on $\g_-$. If in addition the commutator in $\g$ is continuous
in this topology then the triple $(\g,\g_+,\g_-)$ is called
a Manin triple. 
 
To every Lie bialgebra $\a$ one can
associate the corresponding Manin triple \linebreak
$(\g=\a\oplus\a^*,\a,\a^*)$,
where the Lie structure on $\g$ is as above. Conversely, if $(\g,\g_+,\g_-)$
is a Manin triple then $\g_+$ is naturally a Lie bialgebra:
the pairing $\<,\>$ identifies $\g_+^*$ with $\g_-$, which defines a commutator
on $\g_+^*$. This commutator turns out to be dual to a 1-cocycle
(cf. \cite{Dr1}). 

Thus, there is a one-to-one correspondence between Lie bialgebras
and Manin triples. 

Let $(\g,\g_+,\g_-)$ be a Manin triple. 
Let $\{a_i,i\in I\}$ be a basis of $\g_+$, 
and $b^i\in \g_-$ be the linear functions on $\a$ defined by
$b^i(a_j)=\delta_{ij}$. 

\proclaim{Lemma 7.1} Let $M$ be an equicontinuous $\g$-module. Then
for any $v\in M$ and any neighborhood of zero $U\subset M$
 one has $b^iv\in U$  
for all but finitely many $i\in I$.
\endproclaim

\demo{Proof} We assume that $\text{dim}\g_+=\infty$ 
(otherwise there is nothing to prove).

Let $\{i_m\in I:m\ge 1\}$ be any
sequence of distinct elements. The $b^{i_m}\to 0$, $m\to\infty$, so 
$b^{i_m}v\to 0$, $m\to\infty$, for any $v\in M$. This means that
$b^iv\in U$ for almost all $i$. 
$\square$\enddemo

{\it 7.5. Examples of equicontinuous $\g$-modules.}

In this section we will construct examples
of equicontinuous $\g$-modules in the case when $\g$ belongs to a Manin
triple $(\g,\g_+,\g_-)$.

Consider the Verma modules $M_+=\text{Ind}_{\g_+}^\g \bold 1$
$M_-=\text{Ind}_{\g_-}^\g \bold 1$, 
(here $\bold 1$ denotes the trivial 1-dimensional representation).
The modules $M_\pm$ are freely generated over 
$U(\g_\mp)$ by a vector $1_\pm$ such that $\g_\pm 1_\pm=0$, 
and thus are identified (as vector spaces) with $U(\g_\mp)$
via $x1_\pm\to x$. 

Below we show that the module $M_-$ and the module $M_+^*$ dual to $M_+$ in
an appropriate sense are equicontinuous $\g$-modules.

\proclaim{Lemma 7.2} 
The module $M_-$, equipped with the discrete topology,
is an equicontinuous $\g$-module. 
\endproclaim

\demo{Proof} In order to prove the continuity of $\pi_{M_-}(g)$ as a function
on $\g$, we have to check that for
any $v\in M_-$ the space $\g_-v\subset M_-$ is finite dimensional.
One may assume that $v=a_{i_1}a_{i_2}...a_{i_n}1_-$.
We show that $\g_-v$ is finite dimensional by induction in $n$. 
The base of induction is clear since $\g_-v=0$ if $n=0$. 
Now assume that $v=a_jw$, where $w=a_{i_1}...a_{i_{n-1}}1_-$. 
By the induction assumption, we know that $\g_-w$ is finite
dimensional. For any $b\in\g_-$  
we have $bv=ba_jw=[ba_j]w+a_jbw$. For any $j\in I$ we denote
by $W_j\subset \g_-$ space of all $b\in \g_-$ such that
$(1\o b)(\delta(a_j))=0$. It is clear that $W_j$
has finite codimension. For any $b\in W_j$, we have 
$[ba_j]\in\g_-$, since 
$\text{ad}^*b(a_j)=0$ by the definition of $W_j$.
Therefore, for any $b\in W_j$ $bv=[ba_j]w+a_jbw\in \g_-w\oplus a_j\g_-w$.
The latter space is finite dimensional, which implies that $W_jv$ is finite
dimensional. Since $W_j$ has a finite 
codimension in $\g_-$, the space $\g_-v$ is finite-dimensional. 
This implies the continuity
of the homomorphism $\pi_{M_-}:\g\to\End M_-$.
The equicontinuity condition is trivial. 
$\square$\enddemo

Let us now introduce a topology on the space 
$M_+$. This topology comes from the identification of $M_+$ with 
$U(\g_-)$. The space $U(\g_-)$ can be represented as a union  
of $U_n(\g_-)$, $n\ge 0$, where $U_n(\g_-)$ is the set of all elements
of $U(\g_-)$ of degree $\le n$. Furthermore, for any $n\ge 0$, 
we have a linear map $\g_-^{\o n}\to U_n(\g_-)$ given by 
$x_1\o...\o x_n\to x_1...x_n$. This map induces a linear isomorphism
$\xi_n:\oplus_{j=0}^nS^j\g_-\to U_n(\g_-)$, 
where $S^j\g_-$ is the $j$-th symmetric
power of $\g_-$ (as usual we set $\g_-^{\o 0}=
S^0\g_-= k$). Since $S^j\g_-$ has a natural weak topology, coming from its 
embedding to $(\g_+^{\o j})^*$, the isomorphism $\xi_j$ defines a 
topology on $U_n(\g_-)$. Moreover, by the definition, if $m<n$ then
$U_m(\g_-)$ is a closed subspace in $U_n(\g_-)$.
This allows us to
 equip $U(\g_-)$, i.e. $M_+$, with the topology of inductive limit.
By the definition, 
a set $U\subset U(\g_-)$ is open in this topology if and only if 
$U\cap U_n(\g_-)$ is open for all $n$.

\proclaim{Lemma 7.3} Let $g\in\g$. Then $\pi_{M_+}(g)$ is a continuous operator
$M_+\to M_+$. 
\endproclaim

\demo{Proof} Let $g\in\g$. We need to show that for any neighborhood
of the origin $U\subset M_+$ there exists a neighborhood of
the origin $U'\subset M_+$ such that $\pi_{M_+}(g)U'\subset U$. 

Let $U\in U(\g_-)$ be a neighborhood of zero, and $U_n=U\cap U_n(\g_-)$. 
To construct $U'$, we need to construct $U^{\prime}_n=
U^{\prime}\cap U_n(\g_-)$ such that $U^{\prime}_n=U^{\prime}_{n+1}
\cap U_n(\g_-)$. Before giving the construction of $U^{\prime}_n$, 
we make some definitions.   

For any neighborhood $U$ of zero, there exists an increasing
sequence of finite subsets
 $T_n\subset I$, $n\ge 1$, such that for any 
$f\in S^m\g_-$, $m\le n$ satisfying 
the equation $f(a_{i_1},...,a_{i_m})=0$
for any $i_1,...,i_m\in T_n$, one has $\xi_n(f)\in U$.
Fix such a sequence $\{T_n,n\ge 1\}$. 

Let $I$ be as in Section 7.4. For any finite subset
$J\subset I$ denote by $S(J)$ the set of all
$i\in I$ such that there exists $b\in\g_-$ and
$j\in J$ with the property $[bb^i](a_j)\ne 0$.
Since $[bb^i](a_j)=b\o b^i(\delta(a_j))$, the set 
 $S(J)$ is finite. Let the sets $S_n(J)\subset I$ be
defined recursively by $S_0(J)=J$, $S_n(J)=S(S_{n-1}(J))$. 

To construct $U'$, we consider separately the cases $g\in \g_+$ 
and $g\in \g_-$. First consider the case $g\in\g_-$. 

For any elements $x_1,...,x_n\in\g_-$ ($n\ge 1$) consider the element
$X=\sum_{\sigma\in S_n}x_{\sigma(1)}...x_{\sigma(n)}$ in $U_n(\g_-)$,
where $S_n$ is the symmetric group. Consider the
element $gX\in U_{n+1}(\g_-)$. It is easy to see that it is possible
to write $gX$ as a linear
combination of elements of the form
$\sum_{\sigma\in S_m}y_{\sigma(1)}...
y_{\sigma(m)}$, $y_p\in\g_-$, $0\le m\le n+1$, in such a way that 
$y_p$ are iterated commutators of $g$ and $x_1,...,x_n$, and
the number of commutators involved in 
each term $y_p$ does not exceed $n$.

Now we make a crucial observation. 

{\bf Claim.} Let $J\subset I$ be a finite subset.
If for some $m$, $1\le m\le n$, we have $x_m(a_i)=0$, for all 
$i\in S_n(J)$, then 
every monomial $y_1...y_m$ in the symmetrized expression   
of $gX$ contains a factor $y_p$ such that $y_p(a_i)=0$, $i\in J$.

{\it Proof.} Clear.

The construction of $U'$ is as follows. 
For $n\ge 1$, let $U^{\prime}_{n}\subset U_{n}(\g_-)$ be the span
 of all elements $\xi_m(f)$, $0\le m\le n$, where
$f\in S^m\g_-$ are such that 
$f(a_{i_1},...,a_{i_m})=0$ whenever $i_1,...,i_m\in S_n(T_{n+1})$.
Also, set $U^{\prime}_0=0$ (recall that $\{0\}\subset k$ is a neighborhood
of zero since $ k$ is discrete).
Our observation shows that
 for any $X\in U^{\prime}_n$, $gX\in U_{n+1}$, as desired.

Now consider the case $g\in \g_+$. 

Let $R_0(g)\subset I$ be the set
of all $i\in I$ such that $b^i(g)\ne 0$. This is a finite set.
Define inductively the sets $R_n(g)$ by $R_n(g)=S(R_{n-1}(g))$. 

For any finite subsets $K,J\subset I$ 
denote by $P(K,J)$ the set of all $i\in I$ such that 
there exists $j\in J$ and $k\in K$ with $[a_kb^i](a_j)\ne 0$. 
It is clear that if $K,J$ are finite then $P(K,J)$ is finite.
Let $P_n(K,J)$ be defined inductively by $P_n(K,J)=P(K,P_{n-1}(K,J))$.
 
Let $n\ge 1$ be an integer, $X\in U_n(\g_-)$ be as above, and 
$K=R_n(g)$. Consider the vector $gX1_+\in M_+$. Using the relations in $M_+$, 
we can reduce this vector to a linear combination of vectors of the form
$\sum_{\sigma\in S_m}y_{\sigma_1}...y_{\sigma_m}$, $y_p\in\g_-$,
$0\le m\le n+1$, in such a way that $y_p$ are obtained by iterated 
commutation of $g$, $x_1,...,x_n$. As before, it is easy to see 
that the resulting symmetrized expression will contain no more than $n$ 
commutators.

Now let us make a crucial observation. 

{\bf Claim.}
Let $J\subset I$ be any finite subset.
If for some 
$m$, $1\le m\le n$, we have $x_m(a_i)=0$, for all $i\in S_n(P(K,S_n(J))$, then 
every monomial $y_1...y_m$ in the symmetrized expression   
of $gX1_+$ contains a factor $y_p$ such that $y_p(a_i)=0$, $i\in J$.

{\it Proof.} Clear.

The construction of $U'$ is as follows.
For $n\ge 1$, let $U^{\prime}_{n}\subset U_{n}(\g_-)$ be the span
 of all elements $\xi_m(f)$, $f\in S^m\g_-$, $0\le m\le n$, such that 
$f(a_{i_1},...,a_{i_m})=0$ whenever $i_1,...,i_m\in S_n(P(K,S_n(T_{n+1})))$.
Also, set $U^{\prime}_0=0$.
Our observation shows that
 for any $X\in U^{\prime}_n$, $gX\in U_{n+1}$, as desired.
$\square$\enddemo
     
Consider the vector space $M_+^*$ of continuous linear functionals on $M_+$. 
By the definition,
$M_+^*$ is naturally isomorphic to the projective limit of $U_n(\g_-)^*$
as $n\to\infty$. As vector spaces, $U_n(\g_-)^*=(S^j\g_-)^*=S^j\g_+$.
Therefore, it is natural to put the discrete topology on 
$U_n(\g_-)^*$. This equips the module 
$M_+^*$ with a natural structure of a 
complete vector space. It is also equipped with
a filtration by subspaces $(M_+^*)_n=U_{n-1}(\g_-)^\perp$, $n\ge 1$,
such that $M_+=\underleftarrow{\lim}M_+^*/(M_+^*)_n$. 

{\bf Remark.} The
topology of projective limit on $M_+^*$ does not, in general, 
coincide with the weak topology of the dual. In fact, it 
is stronger than the weak topology.

By Lemma 7.3, $M_+^*$ has a natural structure of a $\g$-module.
Namely, the action of $\g$ on $M_+^*$ is defined to be the dual to
the action of $\g$ on $M_+$. 

\proclaim{Lemma 7.4} $M_+^*$ is an equicontinuous $\g$-module.
\endproclaim

\demo{Proof} It is easy to see that $a (M_+^*)_n\subset (M_+^*)_n$,
$a\in\g_+$, and $b (M_+^*)_n\subset (M_+^*)_{n-1}$, $b\in\g_-$.
This means that the operators $\pi_{M_+^*}(g)$ are
continuous for any $g\in\g$, 
and $\pi_{M_+^*}(\g)\subset \End M_+^*$ is an equicontinuous
family of operators. It remains to show that the assignment
$g\to \pi_{M_+^*}(g)$ is continuous for $g\in\g$. Since $\g_+$ is discrete,
it is enough to check this statement for $g\in \g_-$.  

Let $f\in M_+^*$. Let $f_n$ be the reduction of $f$ modulo $(M_+^*)_n$.
We can regard $f$ as an element of $\oplus_{j=0}^nS^j\g_+$. Let us write
$f_n$ in terms of the basis $\{a_i\}$, and let $T_n(f)$ be the set of all
$i\in I$ such that $a_i$ is involved in this expression. 

Let $S_n(J)$ be as in the proof of Lemma 7.3, and $i\in I\setminus 
S_n(T_{n+1}(f))$. 
Then it is easy to see that $b^if\in (M_+^*)_n$. This shows that
for any $n\ge 0$ and any $f\in M_+^*$ $b^if\in (M_+^*)_n$ for almost all 
$i\in I$. 

Thus, $M_+^*$ is an equicontinuous $\g$-module.
$\square$\enddemo

{\bf Remark.} If $\g_+$ is infinite dimensional then $M_+$ is not, in general,
an equicontinuous $\g$-module, since the family of operators
$\{\pi_{M_+}(g),g\in\g_+\}$ may fail to be equicontinuous.

\vskip .05in
7.6. {\it The Casimir element.}

Consider the tensor product 
$\a\o \a^*$. This space can be embedded into
$\End \a$, by $(x\o f)(y)=f(y)x$, 
$x,y\in\a$, $f\in\a^*$. This embedding defines a topology 
on $\a\o \a^*$, obtained by restriction of the weak topology on $\End \a$.
Let $\a\tlo\a^*$ be the completion of $\a\o \a^*$ in this topology.
Since the image of $\a\o \a^*$ is dense in $\End \a$, 
this completion is identified with $\End \a$.

\proclaim{Lemma 7.5} Let $V,W\in \Cal M^e_0$. The map
$\pi_V\o\pi_W: \a\o\a^*\to\End (V\ho W)$ extends to a continuous map
$\a\tlo\a^*\to \End(V\ho W)$.
\endproclaim

\demo{Proof} Let $x\in V\ho W$ be a vector. 
It is easy to see that the map 
$\pi_V\o\pi_W(\cdot)x: \a\o\a^*\to V\ho W$ is continuous.
Since the space $V\ho W$ is complete, this map extends to
a continuous map $\a\tlo\a^*\to V\ho W$. This allows us
to define a linear map
$\pi_V\o\pi_W: \a\tlo\a^*\to \End(V\ho W)$.
We would like to show that this map is continuous. 

Let $x\in V\ho W$ be a vector, and $n\ge 0$ be an integer.
Let $P\subset V\ho W$ be an open subspace,
and $U=\{A\in \End(V\ho W): Ax\in P\}$. Since open sets of this form 
generate the topology on $\End(V\ho W)$, it is enough 
to show that there exists a neighborhood of zero $Y\subset \a\tlo\a^*$
such that $(\pi_V\o\pi_W)(Y)\subset U$, i.e. 
$(\pi_V\o\pi_W)(Y)x\subset P$. 

We can assume that $P=V_1\ho W+W_1\ho V$, where
$V_1,W_1$ are open subspaces of $V,W$. 
By the equicontinuity of 
$\pi_V(g),\pi_W(g)$, $g\in\g$, 
there exist open subspaces $V_2\subset V$, $W_2\subset W$
 such that $\pi_V(\g)V_2\subset V_1$, 
$\pi_W(\g)W_2\subset W_1$. 
Let $y\in V\o W$ be a vector in the usual tensor product of $V$ and $W$
such that $y-x\in V_2\ho W+V\ho W_2$. Then
for any $t\in\a\tlo\a^*$ $(\pi_V\o\pi_W)(t)(y-x)\in P$, so
it is enough to find $Y$ satisfying the condition 
$(\pi_V\o\pi_W)(Y)y\subset P$. 

We have 
$y=\sum_{j=1}^m v_j\o w_j$, $v_j\in V, w_j\in W$. Let $X\subset \a$ be a
finite-dimensional
subspace such that for any $b\in X^\perp\subset \a^*$ $bw_j\in W_1$
for $j=1,...,m$.
Such a subspace exists by Lemma 7.1. The set $Y=\a\tlo X^\perp$ 
(the completion of $\a\o X^\perp$ in $\a\tlo\a^*$)
is open in $\a\tlo\a^*$,
and $(\pi_V\o \pi_W)(Y)y\subset P$, 
as desired. This shows the continuity
of $\pi_V\o\pi_W$ on $\a\tlo\a^*$. 
$\square$\enddemo 

Let $r\in\a\tlo \a^*$ be the vector corresponding
to the identity operator under the identification $\a\tlo \a^*$ with $\End \a$.
Let $r^{op}\in\a^*\tlo \a$ be the element  obtained
from $r$ by permutation of the components.
We define the Casimir element
$\Omega\in \a\tlo \a^*\oplus\a^*\tlo \a$. 
to be the sum $r+r^{op}$. 

It is easy to see that, 
$r=\sum a_i\o b^i$, $r^{op}=\sum b^i\o a_i$,
$\Omega=\sum(a_i\o b^i+b^i\o a_i)$. 
 
Let $V,W$ be equicontinuous $\g$-modules, and denote by
$\pi_V:\g\to \End V$,
$\pi_W:\g\to \End W$ the corresponding linear maps.
Let $\Omega_{VW}=\pi_V\o \pi_W(\Omega)$. This endomorphism of $V\ho W$
is well defined and continuous by Lemma 7.5. Moreover, it is easy to
see that $\Omega_{VW}$ commutes with $\g$, so it is 
an endomorphism of $V\ho W$ as an equicontinuous $\g$-module. 

{\bf Remark.} Although the Casimir operator 
$\Omega=\sum (a_i\o b^i+ b^i\o a_i)$ is defined in the product of
any two equicontinuous $\g$-modules $V\ho W$, the Casimir element
$C=\sum (a_ib^i+b^ia_i)$ in general (for $\dim\a=\infty$)
has no meaning as an operator in an equicontinuous
$\g$-module $V$. 
\vskip .05in

7.7. {\it Drinfeld category.}
Let $\Cal M^e$ denote the category whose objects are equicontinuous
$\g$-modules, and
$\text{Hom}_{\Cal M^e}(U,W)=\text{Hom}_\g(U,W)[[h]]$.
This is an additive category. For brevity we will later
write $\text{Hom}$ for $\text{Hom}_{\Cal M^e}$.
 
Define a structure of a braided monoidal category on $\Cal M^e$ 
analogously to Section 1.4, using an associator $\Phi$ and the functor $\ho$.
As before, we identify $\Cal M^e$ with a strict category and forget about 
positions of brackets.

Let $\gamma$ be the functorial isomorphism defined by
$\gamma_{XY}=\beta_{YX}^{-1}\in \Hom(X\o Y,Y\o X)$, $X,Y\in\Cal M^e$.
It is easy to check that $\gamma$ is a braiding on $\Cal M^e$. 
We will need the braiding $\gamma$ in our construction below.

\vskip .1in
\centerline{\bf 8. The fiber functor.}
\vskip .1in

8.1. {\it The category of complete $ k[[h]]$-modules.}

Let $V$ be a complete vector space over $ k$. Then the space 
$V[[h]]=V\ho k[[h]]$ 
of formal power series in $h$ with coefficients in $V$ 
is also a complete vector space. 
Moreover, $V[[h]]$ has a natural
structure of a topological $ k[[h]]$-module.
We call a topological $ k[[h]]$-module complete if it is isomorphic 
to $V[[h]]$ for some complete $V$. 

Let $\Cal A^c$ be the category of complete $ k[[h]]$-modules,
where morphisms are continuous $ k[[h]]$-linear maps. It is an 
additive category. Define the tensor structure on $\Cal A^c$ as follows.
For $V,W\in \Cal A^c$ define $V\tio W$  to be the quotient
of the completed tensor product $V\ho W$
by the image of the operator $h\o 1-1\o h$.
It is clear that for $V,W\in \Cal A^c$, $V\tio W$ is also
in $\Cal A^c$. 
The category $\Cal A^c$ equipped with the functor
$\tio$ is a symmetric monoidal category.
 
Let CVect be the category of complete vector spaces. We have the functor
of extension of scalars, $V\mapsto V[[h]]$, acting from 
$\text{CVect}$ to $\Cal A^c$. This functor respects 
the tensor product, i.e. $(V\ho W)[[h]]$ is naturally
isomorphic to $V[[h]]\tio W[[h]]$. 

\vskip .05in

8.2. {\it Properties of the Verma modules.}

Let $(\g,\g_+,\g_-)$ be a Manin triple, and $\Cal M^e$ be the
Drinfeld category associated to $\g$. Let $M_+$, $M_-$ be the Verma modules
over $\g$ defined in Section 7.5.

Recall that the modules $M_\pm$ are identified with $U(\g_\mp)$. 
Thus, we can define the maps $i_\pm:M_\pm\to 
M_\pm\o M_\pm$ given by comultiplication in the universal 
enveloping algebras $U(\g_\mp)$. These maps are  
$U(\g)$-intertwiners, since they are
$U(\g_\pm)$-intertwiners and map the vector $1_\pm$ to the
$\g_\mp$-invariant vector $1_\pm\o 1_\pm$. 

Let $M_+^*$ be as in Section 7.5, and $f,g\in M_+^*$.
Consider the linear functional $M_+\to k$ defined by
$v\to (f\o g)(i_+(v))$. It is easy to check that this functional is 
continuous, so it belongs to $M_+^*$. Define the map 
$i_+^*: M_+^*\o M_+^*\to M_+^*$ by
$i_+^*(f\o g)(v)=(f\o g)(i_+(v))$, $v\in M_+$. 
It is clear that $i_+^*$ is continuous, so it extends
to a morphism in $\Cal M^e$: $i_+^*: M_+^*\ho M_+^*\to M_+^*$.

Let $V\in \Cal M$. Consider the space
 $\text{Hom}_\g(M_-,M_+^*\ho V)$,
where $\text{Hom}_\g$ denotes the set of continuous homomorphisms.
Equip this space with the weak topology (see Section 7.1).

\proclaim{Lemma 8.1} 
The complete vector space
 $\text{Hom}_\g(M_-,M_+^*\ho V)$ is isomorphic to $V$.
The isomorphism is given by $f\to (1_+\o 1)(f(1_-))$, 
$f\in\text{Hom}_\g(M_-,M_+^*\ho V)$.
\endproclaim

\demo{Proof} By the Frobenius reciprocity, 
$\Hom_\g(M_-,M_+^*\ho V)$ is isomorphic, as a topological vector space,
 to the space of invariants
$(M_+^*\ho V)^{\g_-}$, via $f\to f(1_-)$. Consider the space
$\Hom_{ k}(M_+,V)$ of continuous homomorphisms from $M_+$ to $V$,
equipped with the weak topology,
and the map
$\phi:(M_+^*\ho V)\to \Hom_{ k}(M_+,V)$, given by 
$\phi(f\o v)(x)=f(x)v$, $u\in M_+^*$, $x\in M_+$, $v\in V$.
It is clear that $\phi$ is injective and continuous. 

{\bf Claim.} The map $\phi$ restricts to an isomorphism
$(M_+^*\ho V)^{\g_-}\to \Hom_{\g_-}(M_+,V)$. 

{\it Proof.} It is clear that $\phi((M_+^*\ho V)^{\g_-})\subset
 \Hom_{\g_-}(M_+,V)$. So it is enough to show
that any continuous
$\g_-$-intertwiner $g:M_+\to V$ is of the form $\phi(g')$, 
$g'\in (M_+^*\ho V)^{\g_-}$, where $g'$ 
continuously depends on $\g$.

Let $X\subset V$ be an open subspace. Then for any 
$\g_-$-intertwiner $g:M_+\to V$ and $n\ge 1$
the image of $g(U_n(\g_-)1_+)$ in $V/V_m$ is finite-dimensional.
This shows that $g=\phi(g')$ for some $g'\in (V\o M_+^*)^{\g_-}$. 
It is clear that $g'$ is continuous in $g$. The claim is proved.

By the Frobenius reciprocity, 
 the space $\Hom_{\g_-}(M_+,V)$ is isomorphic to 
$V$ as a topological vector space, via $f\to f(1_+)$. The lemma is proved.
$\square$\enddemo

\vskip .05in
8.3. {\it The forgetful functor.}

Let $F:\Cal M^e\to\Cal A^c$ be a functor given by $F(V)=\Hom(M_-,M_+^*\ho V)$.
Lemma 8.1 implies that this functor is naturally isomorphic to the 
``forgetful'' functor which associates to every equicontinuous $\g$-module $M$ the
complete $ k[[h]]$-module $M[[h]]$. The isomorphism between these 
two functors is given by $f\to (1_+\o 1)(f(1_-))$, for any $f\in F(M)$. 
Denote this isomorphism by $\tau$. 
\vskip .05in

8.4. {\it Tensor structure on the functor $F$.}
  
 From now on, when no confusion is possible,
 we will denote the tensor product in the categories 
$\Cal M^e$ and $\Cal A^c$ by $\o$, instead of $\ho$ and $\tio$. 

 Define a tensor structure on the functor
$F$ constructed in Section 8.3. 

For any $v\in F(V)$, $w\in F(W)$ define $J_{VW}(v\o w)$ to be the
composition of morphisms:
$$
\gather
M_- @>i_->> M_-\o M_-@>v\o w>> M_+^*\o V\o M_+^*\o W
@>1\o\gamma_{23}\o 1>>
\\ M_+^*\o M_+^*\o V\o W @>i_+^*\o 1\o 1>>M_+^* \o V\o W,
\tag 8.1\endgather
$$
where $\gamma_{23}$ denotes the braiding $\gamma$ acting in the second 
and third components of the tensor product.
That is, 
$$
J_{VW}(v\o w)=(i_+^*\o 1\o 1)\circ (1\o\gamma_{23}\o 1)\circ (v\o w)\circ i_-.
\tag 8.2
$$

\proclaim{Proposition 8.2} The maps $J_{VW}$ are isomorphisms
and define a tensor structure 
on the functor $F$.
\endproclaim

\demo{Proof} 
It is obvious that $J_{VW}$ is an isomorphism since
it is an isomorphism modulo $h$.  

To prove the associativity of $J_{VW}$, we need the following result.

\vskip .03in
{\bf Lemma 8.3} $(i_{-}\o 1)\circ i_{-}=(1\o i_{-})\circ i_{-}$
in $\text{Hom}(M_{-},M_{-}^{\o 3})$; 
$(i_+^*\o 1)\circ i_+^*=(1\o i_+^*)\circ i_+^*$
in $\text{Hom}(M_+^*,(M_+^*)^{\o 3})$. 
\vskip .03in

{\it Proof.} The proof of the first identity
coincides with the proof of Lemma 2.3 in Part I.
To prove the second identity, define $M_+\ho M_+\ho M_+$
to be space of continuous linear functionals on $M_+^*\ho M_+^*\ho M_+^*$.
Since the operators
$\Omega_{ij}\in \End_\g(M_+^*\ho M_+^*\ho M_+^*)$ 
are continuous, one can define the dual operators $\Omega_{ij}^*\in 
\End_\g(M_+\ho M_+\ho M_+)$, and hence the operator $\Phi^*$ dual to $\Phi$.
It is easy to show analogously to the proof of Lemaa 2.3
that $\Phi^*(1_+\o 1_+\o 1_+)=1_+\o 1_+\o 1_+$, which implies 
the second identity of Lemma 8.3.
\vskip .03in

Now we can finish the proof of the proposition. We need to show 
that for any $v\in F(V),w\in F(W),u\in F(U)$
$J_{V\o W,U}\circ (J_{VW}\o 1)(v\o w\o u)=J_{V,W\o U}\circ (1\o J_{WU})
(v\o w\o u)$, i.e. 
$$
\gather
(i_+^*\o 1\o 1\o 1)\circ \gamma_{34,5}\circ (i_+^*\o 1\o 1\o 1\o 1)
\circ \gamma_{23}\circ (v\o w\o u)\circ (i_-\o 1)\circ i_-=\\
(i_+^*\o 1\o 1\o 1)\circ \gamma_{23}\circ (1\o 1\o i_+^*\o 1\o 1)
\circ \gamma_{45}\circ (v\o w\o u)\circ (1\o i_-)\circ i_-
\tag 8.3
\endgather
$$
in $F(V\o W\o U)$, where $\gamma_{23,4}$ 
means the braiding applied to the product of the second and the third
factors and 
to the fourth factor.
Because of Lemma 8.3 and commutation relation of $\gamma_{23,4}$ and 
$i_+^* \o 1\o 1\o 1\o 1$, identity (8.3) is equivalent to the identity
$$
\gather
(i_+^*\o 1\o 1\o 1)\circ (i_+^*\o 1\o 1\o 1\o 1)\circ \gamma_{34,5}
\circ \gamma_{23}=\\
(i_+^*\o 1\o 1\o 1)\circ \gamma_{23}\circ (1\o 1\o i_+^*\o 1\o 1)
\circ \gamma_{45}
\tag 8.4
\endgather
$$
in $\Hom(M_+^*\o V\o M_+^*\o W\o M_+^*\o U,M_+^*\o V\o W\o U)$.

To prove this equality, we observe that the functoriality of the braiding
implies the identity
$$
\gamma_{23}\circ (1\o 1\o i_+^*\o 1\o 1)=(1\o i_+^*\o 1\o 1\o 1)\circ 
\gamma_{2,34}.
\tag 8.5
$$
Using (8.5) and the identity $(i_+^*\o 1)\circ i_+^*=(1\o i_+^*)\circ i_+^*$,
which follows from Lemma 8.3, we reduce (8.4) to 
the identity $\gamma_{34,5}\gamma_{23}=\gamma_{2,34}\gamma_{45}$,
which follows directly from the braiding axioms.
$\square$\enddemo

We will call the functor $F$ equipped with the tensor structure defined above 
{\it the fiber functor}.

\vskip .1in
\centerline{\bf 9. Quantization of Lie bialgebras.}
\vskip .1in

\vskip .05in
9.1. {\it The algebra of endomorphisms of
the fiber functor.}

Let $H=\text{End}(F)$ be the algebra of endomorphisms of the functor $F$,
with a topology defined by the ideal $hH\subset H$.
It is clear that $H$ is a topological algebra over $ k[[h]]$ 
(see Part I, Section 3.1).

Let $H_0$ be the algebra of endomorphisms of the forgetful functor
$\Cal M^e_0\to CVect$. It follows from Lemma 8.1 that the algebra 
$H$ is naturally 
isomorphic to $H_0[[h]]$.

Let $F^2:\Cal M^e\times\Cal M^e\to \Cal A^c$ be the bifunctor defined by 
$F^2(V,W)=F(V)\o F(W)$. Let $H^2=\text{End}(F^2)$.
It is clear that $H^2\supset H\o H$ but $H^2\ne H\o H$
unless $\g$ is finite dimensional.

The algebra $H$ has a natural ``comultiplication'' $\Delta:H\to H^2$
defined by $\Delta(a)_{V,W}(v\o w)=J_{VW}^{-1}a_{V\o W}
J_{VW}(v\o w)$, $a\in H$, $v\in F(V),w\in F(W)$
where $a_V$ denotes the action of $a$ in $F(V)$. We can also define the counit 
on $H$ by $\e(a)=a_{\bold 1}\in k[[h]]$,
where $\bold 1$ is the neutral object.

A topological algebra $A$ over $ k[[h]]$ is said to be a topological bialgebra
if it is equipped with a coproduct $\Delta:A\to A\o A$ 
(where $\o$ is the tensor product in $\Cal A$) and a counit
$\e:A\to k[[h]]$ which are $ k[[h]]$-linear, continuous, and satisfy
the standard axioms of a bialgebra. 

We will need the following statement.

\proclaim{Proposition 9.1} Let $A\subset H$ be a topological 
subalgebra such that $\Delta(A)\subset A\o A$. Then $(A,\Delta,\e)$ is 
a topological bialgebra over $ k[[h]]$.  
\endproclaim

The proof is straightforward.

{\bf Remark.} For infinite-dimensional $\g$, the algebra $H$ equipped 
with the topology defined by the ideal $hH$ 
is not a topological bialgebra 
since $\Delta(H)$ is not a subset of $H\o H$. 

In the following sections we construct a quantum universal enveloping
algebra $U_h(\g_+)$, which is a quantization of the Lie 
bialgebra $\g_+$, in the sense of Drinfeld (see \cite{Dr1} 
and Part I, Section 3.1). Namely, the algebra $U_h(\g_+)$ is obtained
as a subalgebra of $H$ such that $\Delta(A)\subset A\o A$.

\vskip .05in
 
9.2. {\it The algebra $U_h(\g_+)$.} 

Let $x\in F(M_-)$. 
Define the endomorphism $m_+(x)$ of the functor $F$ as follows. 
For any $V\in\Cal M^e$, $v\in F(V)$, define the element $m_+(x)v\in F(V)$
to be the composition
of the following morphisms in $\Cal M^e$: 
$m_+(x)v=(i_+^*\o 1)\circ (1\o v)\circ x$. This defines a linear map
$m_+: F(M_-)\to H$. Denote the image of this map by $U_h(\g_+)$. 

It is easy to see that for any $a\in U(\g_+)$ 
$\tau(m_+(a1_-)v)\equiv a\tau(v)\text{ mod }h$, 
which implies that $m_+$ is an embedding.

\proclaim{Proposition 9.2} $U_h(\g_+)$ is a subalgebra in $H$.
\endproclaim

\demo{Proof} 

 Using Lemma 8.3, for any $x,y\in F(M_-)$, $V\in \Cal M^e$, $v\in F(V)$ we obtain
$$
\gather
m_+(x)m_+(y)v=(i_+^*\o 1)\circ (1\o i_+^*\o 1) \circ (1\o 1\o v)\circ 
(1\o y)\circ x=\\
(i_+^*\o 1)\circ (i_+^*\o 1\o 1) \circ (1\o 1\o v)\circ 
(1\o y)\circ x=\\
(i_+^*\o 1)\circ (1\o v)\circ (i_+^*\o 1) \circ  
(1\o y)\circ x=\\
(i_+^*\o 1)\circ (1\o v)\circ z,\tag 9.1
\endgather
$$
where $z=(i_+^*\o 1) \circ  
(y\o 1)\circ x\in F(M_-)$. 

So by the definition we get $m_+(x)\circ m_+(y)=m_+(z)$.
$\square$\enddemo

Note that the algebra $U_h(\g_+)$ is a deformation of the algebra
$U(\g_+)$. Indeed, we can define a linear isomorphism $\mu:U(\g_+)[[h]]
\to U_h(\g_+)$ by $\mu(a)=m_-(a1_-)$, $a\in U(\g_+)[[h]]$. This isomorphism
has the property $\mu(ab)=\mu(a)\circ \mu(b)\text{ mod }h^2$, 
which follows from the fact that $\Phi\equiv 1\text{ mod }h$, 
but in general $\mu(ab)\ne \mu(a)\circ \mu(b)$.

The subalgebra $U_h(\g_+)$ has a unit which is 
equal to $\mu(1)$, $1\in U(\g_+)$.
To check this, it is enough to observe that 
$\mu(1)$ is invertible and check the identity $\mu(1)^2=\mu(1)$. 

\vskip .05in

9.3. {\it The coproduct on $U_h(\g_+)$.} 

\proclaim{Proposition 9.3} The algebra $U_h(\g_+)$ is closed
under the coproduct $\Delta$, i.e. 
$\Delta(U_h(\g_+))\subset U_h(\g_+)\o U_h(\g_+)$,
and for any $x\in F(M_-)$ one has
$$
\Delta(m_+(x))=(m_+\o m_+)(J^{-1}_{M_-M_-}((1\o i_-)\circ x)).\tag 9.2
$$
\endproclaim

\demo{Proof} Let $x\in F(M_-)$, $V,W\in\Cal M^e$, $v\in V$, $w\in W$.
By the definition of $\Delta$ and $m_+$, the element 
$\Delta(m_+(x))\in H^2$ 
is uniquely determined by the identity
$$
\gather
(i_+^*\o 1\o 1)\circ (1\o i_+^*\o 1\o 1)\circ \gamma_{34}\circ (1\o v\o w)\circ
R(1\o i_-)\circ x=\\
(i_+^*\o 1\o 1)\circ \gamma_{23}\circ \Delta(m_+(x))(v\o w)
\circ i_-\tag 9.3\endgather
$$
in $F(V\o W)$.

The element $X=J_{M_-M_-}^{-1}((1\o i_-)x)\in F(M_-)\o F(M_-)$ 
is, by the definition, uniquely determined by the identity
$$
(1\o i_-)\circ x=(i_+^*\o 1\o 1)\circ \gamma_{23}\circ X\circ i_-\tag 9.4
$$
in $F(M_-\o M_-)$. Therefore, to prove formula (9.2), it is enough to 
prove the equality obtained by substitution
of $(i_+^*\o 1\o i_+^*\o 1)\circ (1\o v\o 1\o w)\circ X$ instead of 
$\Delta(m_+(x))(v\o w)$ in (9.3):
$$
\gather
(i_+^*\o 1\o 1)\circ (1\o i_+^*\o 1\o 1)\circ \gamma_{34}\circ (1\o v\o w)\circ
(1\o i_-)\circ x=\\
(i_+^*\o 1\o 1)\circ \gamma_{23}\circ (i_+^*\o 1\o i_+^*\o 1)\circ 
(1\o v\o 1\o w)\circ X
\circ i_-\tag 9.5\endgather
$$
in $F(V\o W)$.

Using the functoriality of the braiding and Lemma 8.3, we obtain
$$
\gather
(i_+^*\o 1\o 1)\circ \gamma_{23}\circ (i_+^*\o 1\o i_+^*\o 1)\circ 
(1\o v\o 1\o w)=\\
(i_+^*\o 1\o 1)\circ \gamma_{23}\circ (i_+^*\o 1\o i_+^*\o 1)\circ 
\gamma_{23,4}^{-1}\circ 
(1\o 1\o v\o w)\circ \gamma_{23}=\\
(i_+^*\o 1\o 1)\circ (i_+^*\o i_+^*\o 1\o 1)\circ \gamma_{3,45}\circ  
\gamma_{23,4}^{-1}\circ 
(1\o 1\o v\o w)\circ \gamma_{23}=\\
(i_+^*\o 1\o 1)\circ (i_+^*\o i_+^*\o 1\o 1)\circ \gamma_{45}\gamma_{23}^{-1}
\circ (1\o 1\o v\o w)\circ \gamma_{23}=\\
(i_+^*\o 1\o 1)\circ (i_+^*\o 1\o 1\o 1)\circ 
(1\o i_+^*\o 1\o 1\o 1)\circ \gamma_{45}\gamma_{23}^{-1}
\circ (1\o 1\o v\o w)\circ \gamma_{23}\tag 9.6
\endgather
$$
in $\text{Hom}(M_+^*\o M_-\o M_+^*\o M_-,M_+^*\o V\o W)$. 
It is easy to see that $i_+^*\circ \gamma=i_+^*$, so using Lemma 8.3 again,
we get from (9.6):
$$
\gather
(i_+^*\o 1\o 1)\circ \gamma_{23}\circ (i_+^*\o 1\o i_+^*\o 1)\circ 
(1\o v\o 1\o w)=\\
(i_+^*\o 1)\circ 
(1\o i_+^*\o 1\o 1)\circ \gamma_{34}
\circ (1\o v\o w)\circ (i_+^*\o 1\o 1)\circ \gamma_{23}\tag 9.7
\endgather
$$

Substituting (9.7) into the right hand side of (9.5) and using (9.4), we get
$$
\gather
(i_+^*\o 1\o 1)\circ \gamma_{23}\circ (i_+^*\o 1\o i_+^*\o 1)\circ 
(1\o v\o 1\o w)\circ X
\circ i_-=\\
(i_+^*\o 1\o 1)\circ (1\o i_+^*\o 1\o 1)\circ \gamma_{34}
\circ (1\o v\o w)\circ (i_+^*\o 1\o 1)\circ \gamma_{23}
\circ X
\circ i_-=\\
(i_+^*\o 1\o 1)\circ (1\o i_+^*\o 1\o 1)\circ \gamma_{34}
\circ ( 1\o v\o w)\circ (1\o i_-)\circ x
\tag 9.8\endgather
$$
in $F(V\o W)$, which proves (9.2). The proposition is proved.
$\square$\enddemo

\proclaim{Corollary 9.4} The algebra $U_h(\g_+)$,
equipped with the coproduct $\Delta$, is a quantized universal 
enveloping algebra.
\endproclaim

\demo{Proof.} It follows from Lemma 9.1 and 
Propositions 9.2, 9.3 that $U_h(\g_+)$ is a topological bialgebra 
over $ k[[h]]$ isomorphic to $U(\g_-)[[h]]$ as a topological $ k[[h]]$-module,
and such that $U_h(\g_+)/hU_h(\g_+)$ is isomorphic to $U(\g_+)$ as a 
bialgebra. This implies that $U_h(\g)$ has an antipode, because
the antipode exists $\text{mod} h$. Thus, 
$U_h(\g_+)$ is a quantized universal enveloping algebra. 
$\square$\enddemo

9.4. {\it The algebra $U_h(\g_+)$ is a quantization of $\g_+$.}

\proclaim{Proposition 9.5} The algebra 
$U_h(\g_+)$ is a quantization of the Lie bialgebra $\g_+$. 
\endproclaim

\demo{Proof} Let $x\in U_h(\g_+)$ be such that 
there exists $x_0\in\g_+\subset U(\g_+)$ satisfying the 
condition $x\equiv x_0\text{ mod }h$.

It is easy to show that for any $V,W\in M$
$$
\tau_{V\o W}^{-1}\circ J_{VW}\circ (\tau_V\o \tau_W)=
1+hr/2+O(h^2)\tag 9.9
$$ 
in $\End(V\o W)$. From (9.9) and the definition of coproduct, 
analogously to the proof of Proposition 3.6 in Part I, it is easy
to obtain the congruence
$$
h^{-1}(\Delta(x)-\Delta^{op}(x))\equiv \delta(x_0)\text{ mod }h.\tag 9.10
$$
which means that $U_h(\g_+)$ is the quantization of $\g_+$. 
$\square$\enddemo

Thus, we have proved the following theorem, which answers question 1.1
in \cite{Dr3}. 

\proclaim{Theorem 9.6} Let $\a$ be a Lie bialgebra over $ k$.
Then there exists a
quantized universal enveloping algebra $U_h(\a)$ over $ k$ which is a 
quantization of $\a$.
\endproclaim

9.5. {\it The isomorphism between two constructions of the quantization.}

Let us compare the results of the previous sections to
the results of Part I. In Part I, we showed the existence of
quantization for any finite dimensional Lie bialgebra. Let $(\g,\g_+,\g_-)$
be a finite-dimensional Manin triple. Let $U_h(\g_+)$ denote the 
quantization of $\g_+$ constructed in this section, and by $\tilde U_h(\g_+)$
the quantization constructed in Part I.

\proclaim{Proposition 9.7} The quantized universal enveloping algebras
$U_h(\g_+)$, $\tilde U_h(\g_+)$ are isomorphic.
\endproclaim

(Note added on May 19, 2016: Adrien Brochier has discovered that the proof of Proposition 9.7 given below contains an error, 
namely the morphism $\chi$ is defined incorrectly. A corrected proof with the right definition of $\chi$ appears in [Br], Subsection 
2.2). 

\demo{Proof} If $\g$ is finite-dimensional, then $M_+$ is an equicontinuous 
$\g$-module. Let $\tilde F: \Cal M^e\to\Cal A^c$ be
the functor defined by $\tilde F(V)=\Hom(M_+\o M_-,V)$, $V\in\Cal M^e$.
The tensor structure on $\tilde F$ can be defined as in Part I. 

Let $\sigma\in \Hom(\bold 1,M_+^*\o M_+)$ be the canonical element.  
 Consider the morphism $\chi: \tilde F\to F$, defined as follows. 
For any $V\in M$, $v\in \tilde F(V)$, define $\chi_V(v)\in F(V)$ as 
the composition $\chi_V(v)=(1\o v)\circ (\sigma\o 1)$. 
It is obvious that $\chi$ is an isomorphism of additive functors. 

{\bf Claim.} $\chi$ is an isomorphism of tensor functors.

{\it Proof.} The statement is equivalent to the identity
$$
\gather
(1\o v\o w)\circ \beta_{34}\circ (1\o i_+\o i_-)\circ (\sigma\o 1)=\\
(i_+^*\o 1\o 1)\circ \gamma_{23}\circ (1\o v\o 1\o w)\circ 
(\sigma\o 1\o\sigma\o 1)\circ i_-,\tag 9.11\endgather
$$
which should be satisfied in $\Hom(M_-,M_+^*\o V\o W)$
for any $V,W\in\Cal M^e$, $v\in \tilde F(V)$, $w\in\tilde F(W)$. 
Using the identity
$(1\o v\o 1\o w)\circ \gamma_{23}=\gamma_{23,4}\circ (1\o 1\o v\o w)$, 
we reduce (9.11) to the identity
$$
\gather
\beta_{34}\circ (1\o i_+\o i_-)\circ (\sigma\o 1)=\\
(i_+^*\o 1\o 1\o 1\o 1)\circ \gamma_{23,4}\circ 
(\sigma\o 1\o\sigma\o 1)\circ i_-\tag 9.12\endgather
$$
in $\Hom(M_-,M_+^*\o M_+\o M_-\o M_+\o M_-)$.
Moving $\beta_{34}$ from left to right and interchanging $\beta_{34}^{-1}$ 
with $i_+^*\o 1\o 1\o 1\o 1$, so that (9.12) is equivalent to the identity:
$$
\gather
(1\o i_+\o i_-)\circ (\sigma\o 1)=\\
(i_+^*\o 1\o 1\o 1\o 1)\circ \beta_{45}^{-1}\gamma_{23,4}\circ 
(\sigma\o 1\o\sigma\o 1)\circ i_-\tag 9.13\endgather
$$
in $\Hom(M_-,M_+^*\o M_+\o M_+\o M_-\o M_-)$. It is clear that
$\gamma_{1,23}\circ (1\o\sigma)=\sigma\o 1$ in $\Hom(M_-,M_-\o M_+^*\o M_+)$.
Therefore, using the relations
$\gamma_{23,4}\gamma_{3,45}^{-1}=\gamma_{23}\gamma_{45}^{-1}$, and
$\beta\gamma=1$, we reduce (9.13) to
$$
(1\o i_+)\circ \sigma=
(i_+^*\o 1\o 1)\circ \gamma_{23}
\circ (\sigma\o \sigma) \tag 9.14
$$
in $\Hom(\bold 1, M_+^*\o M_+\o M_+)$. Since $i_+^*\circ \gamma=i_+^*$,
we can rewrite (9.14) as
$$
(1\o i_+)\circ \sigma=
(i_+^*\o 1\o 1)\circ \gamma_{12,3}
\circ (\sigma\o \sigma).\tag 9.15
$$
Using the equality $\gamma_{12,3}\circ (\sigma\o 1)=1\o \sigma$,
we reduce (9.15) to
$$
(1\o i_+)\circ \sigma=
(i_+^*\o 1\o 1)\circ (1\o\sigma\o 1)\circ \sigma.\tag 9.16
$$
To prove this equality,
we compute the image of $1\in \bold 1$ under right hand side of (9.16).
In this calculation, we can ignore the action of the associator because 
for any representations $V_1,V_2,V_3$ of $\g$ the associator acts trivially
on the $\g$-invariants in $V_1\o V_2\o V_3$. The calculation 
yields that $1$ goes to $(1\o i_+)(\sigma(1))$, which 
proves  (9.16). The claim is proved.

Let $\Cal M\subset \Cal M^e$ be the full subcategory of 
discrete $\g$-modules, and
$\tilde U_h(\g)=\End(\tilde F|_{\Cal M})$ be the quantization of $\g$ 
constructed in Part I. It is easy to show that the 
homomorphism of topological Hopf algebras $\End \tilde F\to \tilde U_h(\g)$
defined by restriction from $\Cal M^e$ to $\Cal M$
is an isomorphism,
since both algebras are canonically isomorphic to $U(\g)[[h]]$. 
This means that the morphism $\chi$ defined above 
induces an isomorphism of topological Hopf algebras 
$\tilde U_h(\g)\to U_h(\g)$. It is easy to check that this isomorphism maps
$\tilde U_h(\g_+)$ onto $U_h(\g_+)$, which proves the proposition.
$\square$\enddemo

\vskip .1in
\centerline{\bf 10. Universality and functoriality of the quantization of
Lie bialgebras.}
\vskip .1in

10.1. {\it Acyclic functions.} 

Let $V$ be vector space over $k$.
For any integers $m,n\ge 0$, let $H_{mn}=\Hom(V^{\o m},V^{\o n})$
be the space of tensors of rank $m,n$ on $V$.

Let $B=\oplus_{m,n\ge 0}H_{mn}$. We have two binary operations
on $B$: the tensor product and the composition. (If the composition
makes no sense, we set it to zero). 

Let $m_1,...,m_r,n_1,...,n_r$ be nonnegative integers, and
$W=\oplus_{i=1}^rH_{m_in_i}$. Let $p_i: W\to H_{m_in_i}$ be the natural 
projections. 

Let $X$ be a subset of $W$ and $Y_X$ be the space of all
functions from $X$ to $B$. Denote by $A_X$ the smallest subspace in
$Y_X$ closed under composition and tensor product, and
satisfying the following conditions:

(i) $p_i|_X\in A_X$, $i=1,...,r$;

(ii) If $\sigma\in Y_X$ is a permutation operator from $S_p\subset H_{pp}$,
regarded as a constant function on $X$, then
$\sigma\in A_X$.

We call an element of $A_X$ an acyclic function on $X$. 
\vskip .1in

10.2. {\it Universal quantization.}

Let $\g$ be a Lie bialgebra over $k$, and $U_h(\g)$ be its quantization 
constructed in Chapter 9. Recall that as a $k[[h]]$-module, 
$U_h(\g)$ was identified with $U(\g)[[h]]$, which, in turn, we
identify with $S\g[[h]]$ in the standard way. Therefore, the multiplication 
map $\mu: U_h(\g)\o U_h(\g)\to U_h(\g)$ splits in a direct sum of linear maps
$\mu_{mn}^p: S^m\g\o S^n\g\to S^p\g[[h]]$. Similarly, the coproduct
$\Delta: U_h(\g)\to U_h(\g)\o U_h(\g)$
splits in a direct sum of linear maps
$\Delta_p^{mn}: S^p\g\to S^m\g\o S^n\g[[h]]$. All these linear maps
are functions of the commutator $[,]$ and cocommutator $\delta$ 
of the Lie bialgebra $\g$.

Now consider the setting of Section 10.1, with  
$m_1=2,m_2=1,n_1=1,n_2=2$, $W=\Hom(V\o V,V)\oplus \Hom(V,V\o V)$,
$X\subset W$ the set of all pairs $([,],\delta)\in W$ satisfying 
the axioms of a Lie bialgebra. To every Lie bialgebra
$\g\in X$ ($\g=(V,[,],\delta)$), we have associated a quantized
universal enveloping algebra $U_h(\g)$, which is identified 
with $SV[[h]]$ as a $k[[h]]$-module. Thus, we can regard
$\mu_{mn}^p$, $\Delta^{mn}_p$ as functions on $X$ with values
in $H_{m+n,p}[[h]]$, $H_{p,m+n}[[h]]$, respectively. 
 
\proclaim{Theorem 10.1} The coefficients of the $h$-expansion
of $\mu_{mn}^p$, $\Delta_p^{mn}$, are acyclic functions
on $X$.
\endproclaim

{\bf Remark. } Drinfeld calls a quantization having this property
a universal quantization. Thus, the quantization of Lie bialgebras 
constructed in Chapter 9 is universal. 
 
Theorem 10.1 implies functoriality of quantization. Namely, let
LBA denote the category of Lie bialgebras over $k$, and QUEA denote the 
category of quantum universal enveloping algebras over $k[[h]]$. 

\proclaim{Theorem 10.2} There exists a functor $Q: LBA \to QUEA$ 
such that for any $\g\in LBA$ we have $Q(\g)=U_h(\g)$. 
\endproclaim

\demo{Proof} On objects, the functor $Q$ is already defined. Now let us
define it on morphisms. Let $f:\g_1\to \g_2$ be a homomorphism 
of Lie bilagebras. It defines a linear map 
$Q(f): S\g_1[[h]]\to S\g_2[[h]]$.  
By Theorem 10.1, this map defines a homomorphism of Hopf algebras
$U_h(\g_1)\to U_h(\g_2)$. The theorem is proved. $\square$.
\enddemo

\vskip .1in

10.3. {\it Proof of Theorem 10.1.} Let $\g_+$ be a Lie bialgebra, 
$U_h(\g_+)$ be its quantization. We will use the notation of Section 9.

By the definition, $U_h(\g_+)=F(M_-)=\Hom(M_-,M_+^*\ho M_-)$. 
We have the identifications $\xi_\pm: S\g_\mp\to M_\pm$
given by 
$$
\xi_\pm(\Sym(x_1\o...\o x_l))=\Sym(x_1...x_l)1_\pm, x_i\in\g_\mp,\tag 10.1
$$
 and 
$\theta: U_h(\g_+)\to M_-[[h]]$ by
$$
\theta(x)=(1_+\o 1)(x1_-).
$$
They give us identifications
 $\eta=\theta^{-1}\xi_-: S\g_+[[h]]\to U_h(\g_+)$, and
$\xi_+^*: M_+^*\to \hat S\g_+$  (here hat denotes the 
completion by degree). From now on we fixed these identifications and thus 
regard the spaces $M_-,M_+^*,U_h(\g_+)$ as sums of spaces of the form
$S^m\g_+$ or $S^m\g_+[[h]]$. This allows us to make sense of the statement
 that certain maps
between tensor products of these spaces, depending on $[,],\delta$, are 
acyclic functions (on $X$). 
  
To prove the theorem,
we need to show that the maps $\mu$, 
$\Delta$, $S$ are acyclic. 

To implement the proof, we need a few Lemmas.

Let $r\in\g_+\o\g_-$ be the classical $r$-matrix, so 
that $\Omega=r+r^{op}$. 

\proclaim{Lemma 10.3} (i) The map $r: M_-\o M_-\to M_-\o M_-$ is acyclic.

(ii) The maps $r,r^{op}: M_+^*\ho M_-\to M_+^*\ho M_-$ are acyclic.
\endproclaim

\demo{Proof} (i)
For any nonnegative integers $m,n$
consider the mapping 
$\g_+^{\o m}\o \g_+^{\o n}\to S\g_+\o S\g_+$,
given by 
$$
x_1\o...\o x_m\o y_1\o...\o y_n\to 
r(x_1...x_n1_-\o y_1...y_m1_-).\tag 10.2
$$
We need to show that this mapping is an acyclic function. 
We can do this by induction in $N=m+n$. If $N=0$, the operator is zero 
and the statement is clear. Assume the statement is proved for $N=K-1$ and 
let us prove it for $N=K$. Using the relation $[x\o 1+1\o x,r]=\delta(x)$, 
$x\in\g_+$, we can reduce the question to the case $m=K,n=0$. 
In this case, the map is again zero, Q.E.D.

(ii) By the same reasoning as in (i), we get the statement for $r$. 
For $r^{op}$, we reduce the question
to proving that the map $M_+^*\to M_+^*\ho M_-$ given by
$v\to r^{op}(v\o 1_-)$ is acyclic. 

Let $u=\Sym(y_1...y_m)1_+\in M_+$, $y_1,...,y_m\in\g_-$. Let us compute
the expression
$$
X=(u\o 1)(r^{op}(v\o 1_-))\in M_-.
$$ 
We get
$$
X=-(r^{op}(u\o 1))(v\o 1_-)=\sum_i\<L(b^i,y_1,...,y_m)1_+,v\>a_i1_-,\tag 10.3
$$
where $a_i,b^i$ are dual bases of $\g_+,\g_-$, and
$L$ is a polynomial of commutators of $b^i,y_1,...,y_m$ over $\Bbb Q$ which
is symmetric in $b_i,y_1,...,y_m$ and depends only on $m$. 

Using the duality of $\g_+$ and $\g_-$, from (10.3) we get
$$
X=\sum_i\<b^i\o y_1\o...\o y_m,D_L(v)\>a_i1_-\tag 10.4
$$
where $D_L(v)\in S\g_+$ is a linear combination of 
iterated cocommutators applied to $v$. 
This implies that $r^{op}(v\o 1_-)$ is a linear combination of iterated
cocommutators applied to $v$, so the map 
$v\to r^{op}(v\o 1_-)$ is acyclic. 
$\square$\enddemo

For any $x\in M_-$, let $\psi_x: M_-\to M_+^*\ho M_-$ be the morphism such that
\linebreak $(1_+\o 1)(\psi_x1_-)=x$. 

\proclaim{Lemma 10.4} The map $M_-\o M_-\to M_+^*\ho M_-$ defined by
$x\o y\to \psi_xy$ is acyclic.
\endproclaim

\demo{Proof} For any $x\in M_-$, $y\in U(\g_+)$ 
 we have $\psi_xy1_-=\Delta_0(y)\psi_x1_-$, where $\Delta_0$ is
the coproduct in $U(\g_+)$. Since the map $\Delta_0$ is obviously acyclic, 
it suffices to show that the assignment $x\to\psi_x1_-$ is
an acyclic map $M_-\to M_+^*\ho M_-$.

Let $z\in U(\g_-)$. Since the vector $\psi_x1_-$ is $\g_-$-invariant, 
we have 
$$
(z1_+\o 1)(\psi_x1_-)=(1_+\o S_0(z))(\psi_x1_-)=S_0(z)x,
$$ 
where $S_0$ is the antipode of $U(\g_-)$.
Let $z=\Sym(z_1\o...\o z_m), x=\Sym(x_1\o...\o x_n)1_-$, $z_i\in\g_-$,
$x_i\in\g_+$. Computing the product $S_0(z)x$, we see that it is
 a linear combination of products of expressions of the form $x_i$ and
$Z=[\ad^* x_{i_1}...\ad^* x_{i_r}z_j,x_i]_+$, applied to $1_-$, 
where $[za]_+$ denotes the
$\g_+$-component of $[za]$. Using the identity $[za]_+=\<1\o z,\delta(a)\>$,
we can rewrite $Z$ in the form 
$Z=\pm\<1\o z_j,(\ad x_{i_r}...\ad x_{i_1}\o 1)\delta(x_i)\>$. This shows that 
the product $S_0(z)x$, regarded as an element of $S\g_+$, 
 can be represented as a linear combination 
of summands of the form 
$$
\<z_1\o...\o z_m\o 1^{\o s}, x'\>,
$$
where $x'\in S\g_+$ is a polylinear symmetric function of $x_1,...,x_n$,
such that the assignment $\Sym(x_1\o...\o x_n)\to x'$ is acyclic. 
This proves the acyclicity of $\psi_x1_-$. $\square$
\enddemo

Now we can show the acyclicity of the product in $U_h(\g_+)$. According to
Chapter 9, for any $x,y\in M_-$, 
$$
x*y=(1_+\o 1_+\o 1)(\Phi^{-1}(1\o\psi_y)\psi_x1_-)\tag 10.5
$$
Since $\Phi^{-1}$ is a noncommutative formal series of $h\Omega_{12}$,
$h\Omega_{23}$, the acyclicity of the map $x\o y\to x*y$ follows from
Lemmas 10.3 and 10.4. 

To prove the acyclicity of the coproduct $\Delta$, consider the linear operator
$\Cal J\in \End_k(M_-\o M_-)[[h]]$
defined by
$$
\Cal J(x\o y)=(1_+\o 1_+\o 1\o 1)(\Phi_{1,2,34}^{-1}\Phi_{2,3,4}\gamma_{23}
\Phi_{2,3,4}^{-1}\Phi_{1,2,34}(\psi_x1_-\o\psi_y1_-)).\tag 10.6
$$
According to Proposition 9.3, the coproduct on $U_h(\g_+)$, (when $U_h(\g_+)$ 
is identified with $M_-$), is written in the form
$$
\Delta(x)=\Cal J^{-1}i_-(x).\tag 10.7
$$
The map $J$ is acyclic by Lemmas 10.3 and 10.4. Therefore, $\Delta$ is 
acyclic. 
\vskip .1in

10.4. {\it Universal quantization of quasitriangular Lie bialgebras.}

Let $\g$ be a quasitriangular 
Lie bialgebra over $k$, and $U^{qt}_h(\g)$ be its quasitriangular
quantization 
constructed in Chapter 6. As a $k[[h]]$-module, 
$U_h^{qt}(\g)$ was identified with $U(\g)[[h]]$, which we
identify with $S\g[[h]]$ in the standard way. Therefore, the multiplication 
map $\mu: U_h^{qt}(\g)\o U_h^{qt}(\g)\to U_h^{qt}(\g)$ splits in a direct sum of linear maps
$\mu_{mn}^p: S^m\g\o S^n\g\to S^p\g[[h]]$. Similarly, the coproduct
$\Delta: U_h^{qt}(\g)\to U_h^{qt}(\g)\o U_h^{qt}(\g)$
 splits in a direct sum of linear maps
$\Delta_p^{mn}: S^p\g\to S^m\g\o S^n\g[[h]]$, and the quantum $R$-matrix
$R$ splits in ta direct sum of $R_{mn}\in S^m\g\o S^n\g[[h]]$. All these linear maps
are functions of the commutator $[,]$ and the classical $r$-matrix $r$ 
of $\g$.

Now consider the setting of Section 10.1, with
$m_1=2,m_2=0,n_1=1,n_2=2$, $W=\Hom(V\o V,V)\oplus V\o V$,
$X\subset W$ the set of all pairs $([,],r)\in W$ satisfying 
the axioms of a quasitriangular Lie bialgebra. To every 
quasitriangular Lie bialgebra
$\g\in X$ ($\g=(V,[,],r)$), we have associated a quantized
universal enveloping algebra $U_h^{qt}(\g)$, which is identified 
with $SV[[h]]$ as a $k[[h]]$-module. Thus, we can regard
$\mu_{mn}^p$, $\Delta^{mn}_p,R_{mn}$ as functions on $X$ with values
in $H_{m+n,p}[[h]]$, $H_{p,m+n}[[h]]$, $H_{0,m+n}[[h]]$, respectively. 
 
\proclaim{Theorem 10.5} The coefficients of the $h$-expansion
of $\mu_{mn}^p$, $\Delta_p^{mn}$, $R_{mn}$ are acyclic functions
on $X$.
\endproclaim

\demo{Proof} The proof is analogous to the proof of Theorem 10.1
.$\square$\enddemo

Theorem 10.5 implies functoriality of 
quasitriangular quantization. Namely, let
QTLBA denote the category of quasitriangular
Lie bialgebras over $k$, and QTQUEA denote the 
category of quasitriangular quantum universal enveloping algebras. 

\proclaim{Theorem 10.6} There exists a functor $Q^{qt}: QTLBA \to QTQUEA$ 
such that for any $\g\in QTLBA$ we have $Q^{qt}(\g)=U_h^{qt}(\g)$. 
\endproclaim

\demo{Proof} The proof is analogous to the proof of Theorem 10.2
.$\square$\enddemo
\vskip .1in

10.5. {\it Universal quantization of classical $r$-matrices.}

Let $A$ be an associative algebra over $k$ with unit, and $r\in A\o A$ 
be a solution of the classical Yang-Baxter equation.
In Chapter 5, we assigned to $(A,r)$ a solution of the quantum
Yang-Baxter equation $R(r)\in A\o A[[h]]$.

Consider the setting of Section 10.1, with
$m_1=2,m_2=0, m_3=0,n_1=1,n_2=1, n_3=2$, $W=\Hom(V\o V,V)\oplus V\oplus
V\o V$,
$X\subset W$ the set of all triples $(*,1,r)\in W$ such that
$*$ is an associative product, $1$ is a unit, and $r$ satisfies the classical 
Yang-Baxter equation. To every 
$A\in X$ ($A=(V,*,1,r)$), we have associated 
a quantum $R$-matrix $R(r)\in A\o A[[h]]$. Thus, we can regard
$R$ as a function on $X$ with values
in $H_{02}[[h]]$. 
 
\proclaim{Theorem 10.7} The coefficients of the $h$-expansion
of $R$ are acyclic functions
on $X$.
\endproclaim

\demo{Proof} The theorem follows from Theorem 10.5.
$\square$\enddemo

Theorem 10.5 implies functoriality of 
quantization of classical $r$-matrices. 

Namely, let us call a {\it classical
Yang-Baxter algebra} a pair
$(A,r)$, where $A$ is an associative algebra
with unit over $k$, and $r\in A\o A$ satisfies the classical 
Yang-Baxter equation, and a {\it quantum Yang-Baxter algebra}
a pair $(A,R)$, where $A$ is an associative algebra
with unit over $k$, and $R\in A\o A[[h]]$ satisfies the quantum
Yang-Baxter equation.
Let CYBA, QYBA denote the categories of classical, respectively 
quantum, Yang-Baxter algebras.
Morphisms in these categories are algebra homomorphisms 
preserving the unit and $r$ (respectively, $R$). 

\proclaim{Theorem 10.8} There exists a functor $Q^{YB}: CYBA \to QYBA$ 
such that for any $(A,r)\in CYBA$ we have $Q^{YB}(A,r)=(A,R(r))$. 
\endproclaim

\demo{Proof} The proof is analogous to the proof of Theorem 10.2.
$\square$\enddemo
\vskip .1in

10.6. {\it Quantization over a complete local $\Bbb Q$-algebra}

Let $K$ be a local Artinian or pro-Artinian 
commutative algebra over $\Bbb Q$, and $I$ be the maximal ideal in $K$.
Let $k=K/I$ (it is a field of characteristic 0). Let
LBA(K), QTLBA(K), QUEA(K), 
QTQUEA(K) be the categories of Lie bialgebras,
quasitriangular Lie bialgebras, quantum 
universal enveloping algebras, quasitriangular quantum universal
enveloping algebras over $K$ which are
topologically free as $K$-modules and cocommutative modulo $I$. 
Let CYBA(K), QYBA(K) be the categories of 
classical Yang-Baxter algebras, quantum Yang-Baxter algebras,
which are topologically free as $K$-modules and trivial modulo $I$
(i.e. $r=0, R=1$ modulo $I$). 

In Section 10.2, we showed that 
$\mu,\delta$ are series of acyclic functions of 
$[,],h\delta$ with {\it rational} coefficients. Similarly, in Section
10.4, we showed that
 $\mu,\Delta,R$ are series of acyclic functions with rational coefficients
of $[,],hr$, and in Section 10.5 that $R$ is a series of
 acyclic functions with rational
coefficients of $*,1,hr$. Therefore, we can use these formulas to define 
quantization over $K$ (the series will converge in the topology of $K$).
This quantization over $K$ possesses the same functorial properties as
the original quantization over $k[[h]]$.
Thus, we obtain quantization
functors $Q_A: LBA(K)\to QUEA(K)$, $Q_K^{qt}: QTLBA(K)\to QTQUEA(K)$,
$Q_K^{YB}: CYBA(K)\to QYBA(K)$.

\vskip .1in   

\centerline{\bf Appendix: computation of the product in $U_h(\a)$ 
modulo $h^3$.}
\vskip .05in
 
To illustrate the proof of Theorem 10.1, here we 
compute the product in the quantization 
$U_h(\a)$ of a Lie bialgebra $a$ 
modulo $h^3$.  In the text below we always assume
summation over repeated indices.

Let $\{a_i, i\in I\}$ be a basis of $\a$, and $\{b^i\}$ be the topological
basis of $\a^*$ dual to $\{a_i\}$.
Let us write down the commutation relations for the Lie
algebra $\g=\a\oplus\a^*$:
$$
[a_ia_j]=c_{ij}^ka_k, [b^ib^j]=f^{ij}_kb^k,
[a_ib^j]=f_i^{jk}a_k-c_{ik}^jb^k.\tag A1
$$

Let $1_+^*\in M_+^*$ be the functional on $M_+^*$ defined by
$1_+^*(x1_+)=\e(x)$, $x\in U(\a)$.

Let $\{(M_+^*)_n\}$ be the
filtration of $M_+^*$ which was defined in Chapter 7. 

For $x\in U(\a)$, let $\psi_x: M_-\to M_+^*\ho M_-$ be the
$\g$-intertwiner
such that 
$$
\psi_x 1_-\equiv 1_+^*\o x1_- \text{ mod }(M_+^*)_1.
$$
For $x,y\in U(\a)$, we defined the quantized product $z=y\circ x$
to be the element of $U(\a)[[h]]$ such that the operator $\psi_z$ is
the composition 
$$
\gather
M_-@>\psi_x>> M_+^*\ho M_-@>1\o\psi_y>>M_+^*\ho (M_+^*\ho M_-)
@>\Phi^{-1}>>\\
(M_+^*\ho M_+^*)\ho M_-@>i_+^*\o 1>>M_+^*\ho M_-.\tag
A2\endgather
$$

We want to compute the product $a_q\circ a_p$ modulo $h^3$. We fix 
elements $\rho_i\in M_+^*$, $i\in I$, such that 
$\rho_i(1_+)=0$, $\rho_i(b^j1_+)=\delta_i^j$. These elements
are uniquely defined modulo $(M_+^*)_2$.

Let $w^i\in M_-$ be the vectors such that
$$
\psi_{a_p} 1_-\equiv 1_+^*\o a_p1_-+\rho_i\o w^i\text{ mod }
(M_+^*)_2\o M_-.\tag A3
$$
We must have $b^j\psi_{a_p}1_-=0$ for all $j$, so 
$1_+^*\o b^ja_p1_- +b^j\rho_i\o w^i=0$. But $b^j\rho_i(1_+)=
\rho_i(-b^j1_+)=-\delta_i^j$, so we get
$w^i=b^ia_p1_-=-f_p^{ik}a_k1_-$.

Thus we get
$$
\psi_{a_p}1_-\equiv 1_+^*\o a_p1_--f_p^{ik}\rho_i\o a_k1_-
\text{ mod }(M_+^*)_2\ho M_-.\tag A4
$$

Using (A4), we get 
$$
\gather
\psi_{a_q}a_r1_-\equiv (a_r\o 1+1\o a_r)\psi_{a_q} 1_-\equiv \\
1_+^*\o a_r a_q 1_--f_q^{ik}a_r\rho_i\o a_k-
f_q^{ik}\rho_i\o a_ra_k1_-\text{ mod }(M_+^*)_2\ho M_-.\tag A5\endgather
$$

We have
$$a_r\rho_i(b^j1_+)= -\rho_i(a_rb^j1_+)= \rho_i(c_{rk}^jb^k1_+)= c_{ri}^j,
\tag A6
$$
Thus,
substituting (A6) into (A5), we get
$$
\psi_{a_q}a_r1_-\equiv 1_+^*\o a_r a_q 1_--f_q^{ik}c_{ri}^j\rho_j\o a_k1_--
f_q^{ik}\rho_i\o a_ra_k1_-\text{ mod }(M_+^*)_2\ho M_-.\tag A7
$$
In particular, we have
$$
\gather
(1\o \psi_{a_q})\psi_{a_p} 1_-\equiv 
1_+^*\o 1_+^*\o a_qa_p1_-\\
-c_{pi}^jf_q^{ik}1_+^*\o\rho_j\o a_k1_--
f_q^{ik}1_+^*\o\rho_i\o a_pa_k1_--f_p^{ik}\rho_i\o 1_+^*\o
a_ka_q1_-+\\
f_p^{ik}c_{kl}^jf_q^{ls}\rho_i\o\rho_j\o a_s1_-+
f_p^{ik}f_q^{ls}\rho_i\o \rho_l\o a_ka_s1_-\\
\text{ mod }(M_+^*)_2\ho
M_+^*\ho M_-+M_+^*\ho (M_+^*)_2\ho M_-.\tag A8
\endgather
$$

The definition of an associator implies 
$$
\Phi=1+\frac{h^2}{24}[t_{12},t_{23}]+O(h^3).\tag A9
$$
(see \cite{Dr2},\cite{Dr4}).
This means that the part of the $h^2$-coefficient of $\Phi_{V_1V_2V_3}^{-1}$ 
which belongs to $\a^*\o \a^*\o \a$ is $\frac{1}{24}c_{ij}^kb^i\o b^j\o a_k$.

Now let us apply $\Phi^{-1}$ to both sides of (A8). 
We want to compute the answer in the form $1_+^*\o 1_+^*\o u+...$, $u\in
M_-[[h]]$. To do this, we only need to use the last two terms on the
r.h.s. of (A8) and the $\a^*\o\a^*\o \a$-part of the quadratic term of
$\Phi$. The calculation gives
$$
\gather
\Phi^{-1}(1\o \psi_{a_q})\psi_{a_p} 1_-\equiv 1_+^*\o 1_+^*\o u
\text{ mod }(M_+^*)_1\ho
M_+^*\ho M_-+M_+^*\ho (M_+^*)_1\ho M_-,\\ 
u=a_qa_p1_-+\frac{h^2}{24}(f_p^{ik}f_q^{ls}c_{kl}^jc_{ij}^ma_ma_s
+f_p^{in}f_q^{ls}c_{il}^ra_ra_na_s)1_-\tag A10
\endgather
$$
This shows that 
$$
a_q\circ a_p=a_qa_p+\frac{h^2}{24}(f_p^{ik}f_q^{ls}c_{kl}^jc_{ij}^ma_ma_s
+f_p^{in}f_q^{ls}c_{il}^ra_ra_na_s)+O(h^3).\tag A11
$$
This formula is analogous to the formula  
deduced by Drinfeld \cite{Dr3} (equation 1.1).

It is easy to see that this formula contains only acyclic
monomials. Therefore, this formula is universal. 

 \end

\end